\documentclass[journal]{IEEEtran}

\usepackage{amsmath,amsfonts,amssymb,amsthm,braket}

\newtheorem{theorem}{Theorem}
\newtheorem{lemma}[theorem]{Lemma}
\newtheorem{corollary}[theorem]{Corollary}
\newtheorem{remark}[theorem]{Remark}
\theoremstyle{definition}
\newtheorem{definition}{Definition}

\usepackage{setspace}
\usepackage{algorithm,algpseudocode}
\floatname{algorithm}{\color{black} Algorithm}
\newcommand{\AlgComment}[1]{\textit{\textcolor{gray}{// #1}}}

\usepackage[shortlabels]{enumitem}
\setlist{itemsep=0.05\baselineskip, 
    itemindent=2\parindent, 
    labelsep=0.5em,
    partopsep=0pt, 
    align=right, nosep,
    parsep=0ex, 
    topsep=0.1\baselineskip, 
    leftmargin=0pt, 
    listparindent=\parindent,
}

\usepackage{array}
\usepackage[caption=false,font=normalsize,labelfont=sf,textfont=sf]{subfig}
\usepackage{textcomp}
\usepackage{stfloats}
\usepackage{url}
\usepackage{verbatim}
\usepackage{graphicx}

\usepackage[hidelinks,colorlinks=true,linkcolor=blue,citecolor=blue]{hyperref}
\usepackage{cite}

\usepackage{booktabs,makecell,multirow,array,colortbl,dcolumn,stfloats,diagbox,varwidth}
\usepackage{xspace,xstring,footmisc}
\usepackage[svgnames,dvipsnames,table]{xcolor}
\definecolor{Platinum}{RGB}{229,228,226}

\newcommand{\Rmnum}[1]{\uppercase\expandafter{\romannumeral #1}}

\begin{document}
\title{Tensor Completion via Monotone Inclusion: Generalized Low-Rank Priors Meet Deep Denoisers}

\author{
Peng Chen\textsuperscript{*},
Deliang Wei\textsuperscript{*},
Jiale Yao,
and Fang Li

\thanks{\textsuperscript{*} These authors contributed equally to this work.}
\thanks{Peng Chen, Jiale Yao, and Fang Li are with the School of Mathematical Sciences, East China Normal University, Shanghai 200241, China (e-mail: pengchen2209@stu.ecnu.edu.cn; 52285500019@stu.ecnu.edu.cn; fli@math.ecnu.edu.cn).}
\thanks{Deliang Wei is with the Data Science and AI institute, Johns Hopkins University, Baltimore, MD 21218, America (e-mail: dwei12@jh.edu).}
}

\markboth{Journal of \LaTeX\ Class Files,~Vol.~18, No.~9, September~2020}%
{Tensor Completion via Monotone Inclusion: Generalized Low-Rank Priors Meet Deep Denoisers}

\maketitle

\begin{abstract}
    Missing entries in multi-dimensional data pose significant challenges for downstream analysis across diverse real-world applications. These data are naturally represented as tensors, and recent completion methods integrating global low-rank priors with plug-and-play denoisers have demonstrated strong empirical performance. However, these approaches often rely on empirical convergence alone or unrealistic assumptions, such as deep denoisers acting as proximal operators of implicit regularizers, which generally does not hold. To address these limitations, we propose a novel tensor completion framework grounded in the monotone inclusion paradigm. Within this framework, deep denoisers are treated as general operators that require far fewer restrictions than in classical optimization-based formulations. To better capture holistic structure, we further incorporate generalized low-rank priors with weakly convex penalties. Building upon the Davis–Yin splitting scheme, we develop the GTCTV-DPC algorithm and rigorously establish its global convergence. Extensive experiments demonstrate that GTCTV-DPC consistently outperforms existing methods in both quantitative metrics and visual quality, particularly at low sampling rates. For instance, at a sampling rate of 0.05 for multi-dimensional image completion, GTCTV-DPC achieves an average mean peak-signal-to-noise ratio (MPSNR) that surpasses the second-best method by 0.717 dB, and 0.649 dB for multi-spectral images, and color videos, respectively.
\end{abstract}


\section{Introduction} \label{introduction}

\IEEEPARstart{T}{he} presence of missing entries in data, often resulting from sensor malfunctions, occlusions, or transmission errors, poses a persistent challenge in data analysis~\cite{LIU2025130610}. Tensors, as multi-dimensional arrays, provide a versatile framework for modeling diverse real-world datasets, such as color images~\cite{7859390}, multi-/hyper-spectral images (MSI/HSI)~\cite{SU2025109903}, color videos~\cite{10948572}, and spatio-temporal traffic data~\cite{10595845}. Consequently, tensor completion has become a pivotal research area, attracting significant attention within the scientific community~\cite{https://doi.org/10.1002/nla.2299, 9877877, 10496551}. In this work, we address tensor completion from highly incomplete observations by developing a unified framework grounded in the \textit{monotone inclusion} paradigm, which integrates generalized low-rank priors with deep denoising priors.

The general problem of missing entries in tensor data can be formulated as $\mathcal{Y} = \mathcal{P}_{\Omega} (\mathcal{X})$, where $\mathcal{Y} \in \mathbb{R}^{n_1 \times n_2 \times \cdots \times n_N}$ denotes the observed data with missing entries, $\mathcal{X}$ is the underlying complete tensor, $\Omega$ is the index set of observed entries, and $\mathcal{P}_{\Omega}$ is the projection operator that preserves entries in $\Omega$ and sets all others to zero:
\begin{equation*}
[\mathcal{P}_{\Omega} (\mathcal{X})]_{i_1 i_2 \dots i_N} = \begin{cases}
    \mathcal{X}_{i_1 i_2 \dots i_N} & \text{if } (i_1, i_2, \dots, i_N) \in \Omega, \\
    0 & \text{otherwise}.
\end{cases}
\end{equation*}
To recover $\mathcal{X}$ from $\mathcal{Y}$, a constrained optimization framework is commonly employed:
\begin{equation*}
    \min_{\mathcal{X}} R (\mathcal{X}) \quad \text{s.t.} \quad \mathcal{P}_{\Omega} (\mathcal{X}) = \mathcal{P}_{\Omega} (\mathcal{Y}).
\end{equation*}
Here, $R(\cdot)$ is a regularizer that encodes the intrinsic structural properties of $\mathcal{X}$. In the following, we review related works in two principal directions for designing $R$.

\begin{figure*}
    \centering
    \includegraphics[width=.95\linewidth]{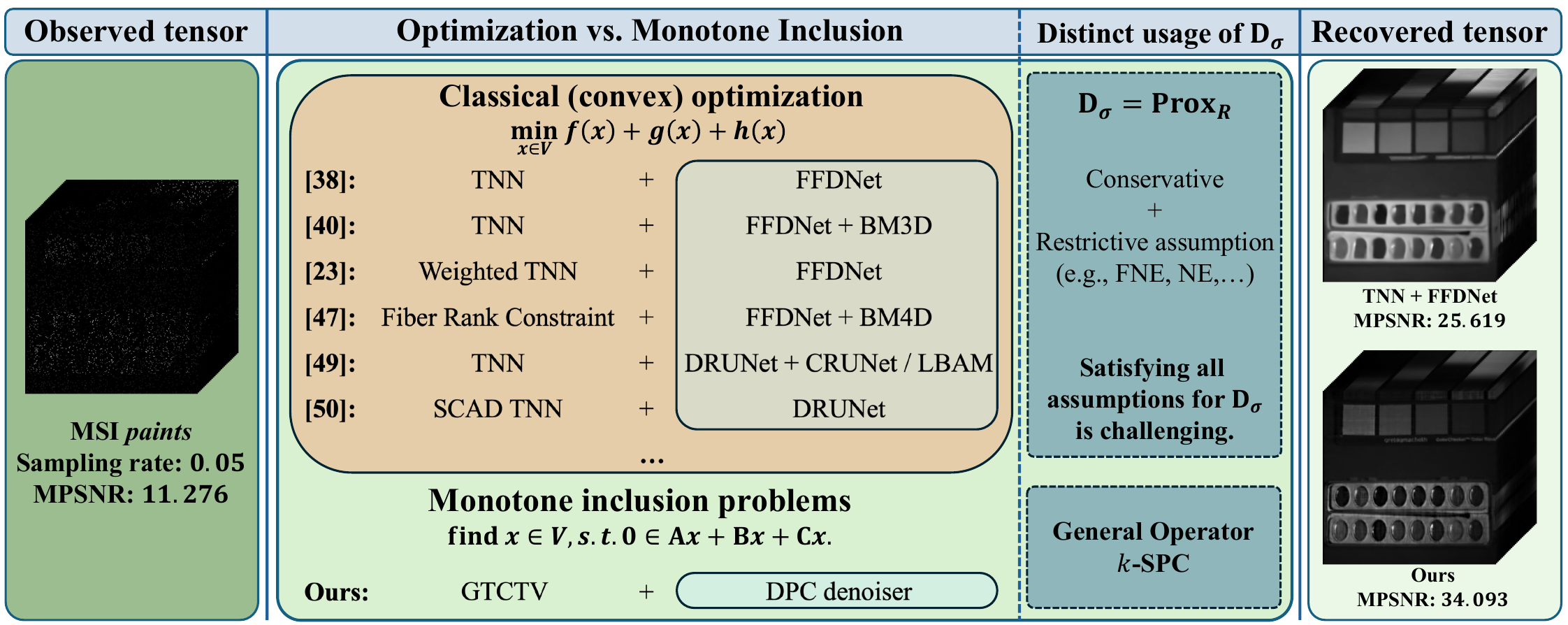}
    \caption{Schematic comparison of recent multi-prior (optimization-based) methods and the proposed monotone-inclusion-based approach. Unlike the former, which treat the denoiser $\operatorname{D}\sigma$ as a proximal mapping under restrictive assumptions, our framework regards deep priors as general operators, imposing much weaker constraints on $\operatorname{D}\sigma$. This relaxation enables more effective denoising~\cite{AHQS} and ultimately yields superior restoration quality.}
    \label{fig:VennMIP}
\end{figure*}

\subsection{Related works} 
\subsubsection{Low-rank based methods}
Low-rankness has emerged as a dominant prior for designing $R$, as it facilitates the extraction of essential structures from high-dimensional data~\cite{11080069}. One prominent approach involves tensor decomposition, such as CANDECOMP/PARAFAC (CP) decomposition~\cite{Kiers2000cp}, Tucker decomposition~\cite{Tucker_1966}, and tensor train decomposition~\cite{Oseledets2011tensortrain}, to represent low-rank properties, with regularization applied to the decomposed components~\cite{TAN201315, 7502115, CHEN201859, CHEN201966, fan2022multimode, Yun2024Inverse, 10530915}. 

Alternatively, tensor nuclear norm (TNN)-based methods, which serve as surrogates for tensor rank, have been shown to outperform decomposition-based approaches~\cite{9115254, 10522984, liu2025tensor}. For multi-dimensional image (MDI) completion, many methods rely on the tubal nuclear norm or its variants derived from tensor singular value decomposition (t-SVD)~\cite{KILMER2011641}. For instance, Jiang et al.~\cite{9115254} employ a framelet transform-based t-SVD and propose a related TNN for visual tensor completion. Additionally, Wang et al.~\cite{10078018} extend this concept by applying TNN to tensor gradients, yielding a tensor correlated total variation (t-CTV) regularizer that jointly exploits low-rankness and smoothness, thereby enhancing completion efficiency. From a functional perspective, Wang et al.~\cite{Wang2025functionalTSVD} generalize t-SVD to functional transforms, resulting in a functional TNN that captures both global low-rank structure and local smoothness. In traffic data completion, Tucker rank approximations via Tucker decomposition remain popular~\cite{doi:10.1137/07070111X, CHEN2020102673, CHEN2024104657}. For instance, Nie et al.~\cite{NIE2022103737} introduce a truncated Schatten $p$-norm (TSpN) to mitigate over-shrinkage of TNN, paired with an adaptive truncation rate decay strategy to handle varying missing rates. 

Recently, deep learning, both supervised and unsupervised, has been leveraged to learn latent low-rank representations through neural architectures~\cite{10510501, 10536168, 10758344}. However, supervised methods require intricate training and often lack generalization, performing poorly on diverse datasets without retraining~\cite{9543471, ZHANG2024102038}. In unsupervised methods, Luo et al.~\cite{10354352} propose a low-rank tensor function representation (LRTFR) parameterized by multi-layer perceptrons based on Tucker decomposition, while Su et al.~\cite{SU2025109903} develop a deep fully-connected tensor network decomposition to capture fine details. Li et al.~\cite{Li_Zhang_Luo_Meng_2025} introduce a simplified deep rank-one tensor functional factorization (DRO-TFF), although its self-supervised performance is limited.

\subsubsection{Multi-prior methods} \label{subsec:related_works}
Beyond global low-rank priors, recent studies integrate additional regularizers to enhance local consistency, often leveraging pretrained deep denoisers in a plug-and-play (PnP) fashion.
PnP frameworks employ a denoiser $\operatorname{D}_\sigma$ to characterize local texture and nonlocal dependencies across slices. While such multi-prior methods improve recovery accuracy, they often lack rigorous theoretical convergence guarantees and rely on empirical validation. 

For instance, Zhao et al.~\cite{zhao2020deep} combine FFDNet~\cite{zhang2018ffdnet} with TNN for MDI completion, achieving significant improvements through the learned prior but providing only empirical convergence. Zhao et al.~\cite{zhao2021tensor} extend this by incorporating both FFDNet and BM3D~\cite{4271520} denoisers within the alternating direction method of multipliers (ADMM)~\cite{boyd2011distributed} framework and prove convergence under the assumption that $\operatorname{D}_\sigma$ is the proximal operator of a Kurdyka--\L{}ojasiewicz (KL) regularizer~\cite{attouch2013convergence}, i.e., $\operatorname{D}_\sigma = \operatorname{Prox}_{R}$. However, recent works~\cite{hurault2022proximal, wei2025learning} demonstrate that deep denoisers are generally not conservative and thus cannot be represented as proximal mappings. Moreover, enforcing denoisers to be proximal often imposes restrictive Lipschitz assumptions, such as residual or firm nonexpansiveness~\cite{hurault2022proximal, pesquet2021learning}, which can compromise denoising performance. Additionally, verifying the KL property of the underlying prior $R$ is challenging in practice. Liu et al.~\cite{liu2025tensor} pair FFDNet with a weighted TNN for anomaly detection in remote sensing images via the ADMM framework, but the convergence analysis still assumes denoisers are proximal operators of implicit regularizers. 

Liang et al.~\cite{liang2024fixed} integrate FFDNet and BM4D~\cite{9165941} denoisers with a fiber-rank constraint and establish fixed-point convergence for a multi-block PnP-ADMM algorithm under bounded-denoiser assumptions. Hu et al.\cite{HU2025112612} integrate two pretrained deep neural networks, one for image denoising and another for completion, with TNN for noisy tensor completion, assuming the networks are nonexpansive and proving fixed-point convergence of their algorithm. Nevertheless, this fixed point may not align with the solution to any optimization objective, thereby limiting interpretability. For traffic data imputation, Chen et al.~\cite{10756233} incorporate a nonconvex tensor low-rank prior and the deep PnP denoiser DRUNet~\cite{zhang2021plug}, providing only empirical convergence analysis for the multi-prior completion method.

\subsection{Contributions}
The synergy between global low-rank priors and deep denoisers has shown considerable potential in tensor completion tasks. However, the limitations highlighted above reveal a gap: low-rank priors effectively capture global structure but often overlook fine local details, whereas optimization-based integrations of deep priors typically assume that $\operatorname{D}_\sigma$ acts as a proximal operator of an implicit regularizer. Such an assumption typically requires the denoiser to be conservative and to satisfy restrictive properties (e.g., firmly nonexpansive), conditions that are difficult to verify in practice and often detrimental to denoising performance. 

This motivates a paragidm shift from \textit{classical optimization} formulations to a broader \textit{monotone inclusion} framework, in which priors are treated as general operators rather than as proximal maps. Within this framework, we can couple weakly convex low-rank regularizers with deep denoisers that satisfy more natural operator properties (e.g., pseudo-contractive), while still obtaining rigorous convergence guarantees. Fig.~\ref{fig:VennMIP} illustrates the distinction between existing multi-prior methods and our proposed approach. Our key contributions are summarized as follows:
\begin{itemize}
    \item \textbf{Tensor completion via monotone inclusion.} We propose a novel tensor completion model formulated within the monotone inclusion framework. The model integrates a generalized tensor correlated total variation (GTCTV) prior to capture global structural dependencies, while incorporating deep pseudo-contractive (DPC) denoisers to preserve fine local details, thereby enhancing overall reconstruction accuracy.  
    \item \textbf{Rigorous global convergence analysis.} Building upon the Davis–Yin splitting (DYS) scheme, we develop the GTCTV-DPC algorithm in Algorithm~\ref{alg:main} and provide a theoretical proof of its global convergence in Corollary~\ref{coro:convergence}.
    \item \textbf{Comprehensive empirical validation.} Extensive experiments demonstrate that GTCTV-DPC consistently outperforms existing methods in both quantitative metrics and visual quality, particularly at low sampling rates. For instance, at a sampling rate of 0.05 for MDI completion, GTCTV-DPC achieves an average MPSNR that surpasses the second-best method by 0.717 dB, and 0.649 dB for MSIs, and color videos, respectively.
\end{itemize}

The remainder of this paper is arranged as follows. Section~\ref{sec:notation} introduces notations and preliminaries. Section~\ref{sec:tc_mip} presents the proposed tensor completion method within the monotone inclusion framework, and provides a rigorous convergence analysis. Section~\ref{sec:experiment} reports numerical experiments validating the proposed method. Finally, section~\ref{sec:conclusion} provides a brief conclusion. 

\section{Notations and preliminaries} \label{sec:notation}
In this section, we introduce the tensor notations used throughout the paper, recall the functional-analytic concepts relevant to our analysis, and give a brief account of the monotone inclusion problem.
\begin{table}[htbp]
    \centering
    \caption{Notations used in this paper.} \label{table:notations}
    \footnotesize{
    \begin{tabular}{c|c}
        \specialrule{.1em}{.05em}{.05em}
        \specialrule{.1em}{.05em}{.5em}
            Notation & Description \\
            \specialrule{.1em}{.5em}{.25em}
            $a, \, \mathbf{a}, \, \mathbf{A}, \, \mathcal{A}$ & Scalar, vector, matrix, tensor. \\
            \specialrule{0em}{.25em}{.25em}
            $\operatorname{diag} \left( \mathbf{a} \right)$ & 
            \makecell{
            The $n$-th order diagonal matrix with $\mathbf{a} \in \mathbb{R}^n$, \\
            where the $(i, \, i)$-th element is $a_{i}$.
            }\\
            \specialrule{0em}{.25em}{.25em}
            $\mathcal{A}_{i_1 i_2 \dots i_N}$ & 
            \makecell{
            The $(i_1, \, i_2, \, \dots, \, i_N)$-th element of \\
            $\mathcal{A} \in \mathbb{C}^{n_1 \times n_2 \times \cdots \times n_N}$.
            }\\
            \specialrule{0em}{.25em}{.25em}
            $\overline{\mathcal{A}_{i_1 i_2 \dots i_N}}$ &
            The complex conjugate of $\mathcal{A}_{i_1 i_2 \dots i_N}$. \\         
            \specialrule{0em}{.25em}{.25em}
            $\mathcal{A}^{(i_3 \dots i_N)}$ & 
            \makecell{
            The $(i_3, \dots, i_N)$-th face slice of $\mathcal{A}$. \\
            $\mathcal{A}^{(i_3 \dots i_N)} := \mathcal{A} \left( :, :, i_3, \, \dots, \, i_N \right)$.
            }\\
            \specialrule{0em}{.25em}{.25em}
            $\mathcal{\langle \mathcal{A}, \mathcal{B} \rangle}$ & 
            \makecell{
            The inner product of $\mathcal{A}$ and $\mathcal{B}$. \\
            $\mathcal{\langle \mathcal{A}, \mathcal{B} \rangle} := \sum_{i_1 i_2 \dots i_N} \overline{\mathcal{A}_{i_1 i_2 \dots i_N}} \mathcal{B}_{i_1 i_2 \dots i_N}$.
            }\\
            \specialrule{0em}{.25em}{.25em}
            $\left\| \mathcal{A} \right\|_F$ & 
            The Frobenious norm of $\mathcal{A}$. $\left\| \mathcal{A} \right\|_F := \sqrt{\mathcal{\langle \mathcal{A}, \mathcal{A} \rangle}}$. \\
            \specialrule{0em}{.25em}{.25em}
            $\mathcal{A} \Delta \mathcal{B}$ & 
            \makecell{Face-wise product of \\
            $\mathcal{A} \in \mathbb{C}^{n_1 \times m \times \cdots \times n_N}$ and $\mathcal{B} \in \mathbb{C}^{m \times n_2 \times \cdots \times n_N}$. \\
            $\mathcal{C} = \mathcal{A} \Delta \mathcal{B} \iff \mathcal{C}^{(i_3 \dots i_N)} = \mathcal{A}^{(i_3 \dots i_N)} \mathcal{B}^{(i_3 \dots i_N)}$. }  \\
            \specialrule{0em}{.25em}{.25em}
            $\operatorname{unfold}_d \left( \cdot \right)$ & 
            \makecell{The unfolding operator along the $d$-th mode. \\
            $\operatorname{unfold}_d: \mathbb{C}^{n_1 \times n_2 \times \cdots \times n_N} \rightarrow \mathbb{C}^{n_d \times \left(\prod_{i \neq d} n_i \right)}$.} \\
            \specialrule{0em}{.25em}{.25em}
            $\operatorname{fold}_d \left( \cdot \right)$ & 
            The inverse operator of $\operatorname{unfold}_d \left( \cdot \right)$. \\
            \specialrule{0em}{.25em}{.25em}
            $\mathbf{A}_{(d)}$ & 
            \makecell{
            The unfolding result of $\mathcal{A}$ along the $d$-th mode. \\
            $\mathbf{A}_{(d)} = \operatorname{unfold}_d (\mathcal{A})$.
            }\\
            \specialrule{0em}{.25em}{.25em}
            $\times_d$ & 
            \makecell{
            Tensor-matrix product along the $d$-th mode. \\
            $\mathcal{A} \times_d \mathbf{F} := \operatorname{fold}_d \left( \mathbf{F} \mathbf{A}_{(d)} \right)$.
            }\\
            \specialrule{0em}{.25em}{.25em}
            $\nabla_d$ & 
            \makecell{The gradient operator along the $d$-th mode. \\
            $\nabla_d \left( \mathcal{A} \right) := \mathcal{A} \times_d \mathbf{D}_{n_d}$, \\
            $\mathbf{D}_{n_d}$ is a row circulant matrix of $\left( -1, 1, 0, \dots, 0 \right)$.}  \\
        \specialrule{.1em}{.25em}{.05em}
        \specialrule{.1em}{.05em}{.05em}
    \end{tabular}}
\end{table}

\subsection{Notations}
In this paper, we adopt the tensor notations from~\cite{10078018, 9730793, doi:10.1073/pnas.2015851118} and focus on order-$N$ tensors with $N \geq 3$. The primary tensor space under consideration is denoted by $\mathbb{V} = \mathbb{R}^{n_1 \times n_2 \times \cdots \times n_N}$, with $\mathcal{O} \in \mathbb{V}$ denoting the zero tensor. Additional notations are summarized in Table~\ref{table:notations}.

Let $\mathfrak{L}$ be an invertible linear transform~\cite{9730793} associated with transform matrices $\{\mathbf{U}_{n_3}, \dots, \mathbf{U}_{n_N}\}$. The transformed representation of a tensor $\mathcal{A}$ is 
\begin{equation*}
    \mathfrak{L}(\mathcal{A}) := \mathcal{A}_{\mathfrak{L}} = \mathcal{A} \times_3 \mathbf{U}_{n_3} \times_4 \cdots \times_N \mathbf{U}_{n_N}, 
\end{equation*}
where each $\mathbf{U}_{n_i}$ is an $n_i \times n_i$ transform matrix satisfying the existence of $l_{n_i} > 0$ and an $n_i \times n_i$ unitary matrix $\mathbf{W}_{n_i}$ such that $\mathbf{U}_{n_i} = l_{n_i} \mathbf{W}_{n_i}$, for $i = 3, \dots, N$. For example, if $\mathbf{U}_{n_i}$ is the unnormalized discrete Fourier transform (DFT), then $l_{n_i} = \sqrt{n_i}$; if it is the discrete cosine transform (DCT), then $l_{n_i} = 1$~\cite{10078018}. The inverse operation is $\mathfrak{L}^{-1}(\mathcal{A}) := \mathcal{A} \times_3 \mathbf{U}^{-1}_{n_3} \times_4 \cdots \times_N \mathbf{U}^{-1}_{n_N}$, satisfying $\mathfrak{L}^{-1}(\mathfrak{L}(\mathcal{A})) = \mathcal{A}$. For a given invertible linear transform $\mathfrak{L}$, we denote $l = \prod_{i=3}^N l_{n_i}$ as the composite scale factor. For instance, when all transform matrices are the unnormalized DFT, $l = \sqrt{\prod_{i=3}^N n_i}$; when all are DCT, $l = 1$. We next recall several standard definitions in the $\mathfrak{L}$-based algebra~\cite{9730793}.

\begin{definition}[tensor-tensor product~\cite{9730793}]
    For tensors $\mathcal{A} \in \mathbb{R}^{n_1 \times m \times n_3 \times \cdots \times n_N}$ and $\mathcal{B} \in \mathbb{R}^{m \times n_2 \times n_3 \times \cdots \times n_N}$, the transform $\mathfrak{L}$-based tensor-tensor product is given by $\mathcal{A} *_{\mathfrak{L}} \mathcal{B} := \mathfrak{L}^{-1}\left(\mathfrak{L}(\mathcal{A}) \, \Delta \, \mathfrak{L}(\mathcal{B})\right)$.
\end{definition}

For the tensor $\mathcal{A} \in \mathbb{C}^{n_1 \times n_2 \times \cdots \times n_N}$, its \textit{conjugate transpose} $\mathcal{A}^{*} \in \mathbb{C}^{n_2 \times n_1 \times \cdots \times n_N}$ satisfies $[\mathfrak{L} (\mathcal{A}^{*})]^{(i_3 \dots i_N)} = [(\mathcal{A}_{\mathfrak{L}})^{*}]^{(i_3 \dots i_N)}$ for all face slices. A tensor $\mathcal{I} \in \mathbb{R}^{n \times n \times \cdots \times n_N}$ is an \textit{identity tensor} if it satisfies $[\mathcal{I}_{\mathfrak{L}}]^{(i_3 \dots i_N)} = \mathbf{I}_n$ for all face slices, where $\mathbf{I}_n$ is the $n \times n$ identity matrix. A tensor $\mathcal{U} \in \mathbb{C}^{n \times n \times \cdots \times n_N}$ is \textit{orthogonal} if $\mathcal{U}^{*} *_{\mathfrak{L}} \mathcal{U} = \mathcal{U} *_{\mathfrak{L}} \mathcal{U}^{*} = \mathcal{I}$. A tensor $\mathcal{A} \in \mathbb{V}$ is \textit{f-diagonal} if each face slice $\mathcal{A}^{(i_3 \dots i_N)}$ is diagonal. 

\begin{theorem}[t-SVD~\cite{9730793}]
    For any tensor $\mathcal{A} \in \mathbb{V}$, it can be decomposed as $\mathcal{A} = \mathcal{U} *_{\mathfrak{L}} \mathcal{S} *_{\mathfrak{L}} \mathcal{V}^{*}$, where $\mathcal{U} \in \mathbb{R}^{n_1 \times n_1 \times \cdots \times n_N}$ and $\mathcal{V} \in \mathbb{R}^{n_2 \times n_2 \times \cdots \times n_N}$ are orthogonal, and $\mathcal{S} \in \mathbb{V}$ is a f-diagonal tensor.
\end{theorem}

The above $\mathfrak{L}$-based t-SVD can be realized by applying the SVD to each face slice of $\mathcal{A}_{\mathfrak{L}}$ in the transform domain, and then mapping the factors back to the original domain via $\mathfrak{L}^{-1}$.

\subsection{Relevant concepts of functional analysis} \label{subsec:functional_analysis}
In this paper, we mainly adopt the relevant concepts of functional analysis as presented in~\cite{bauschke2017correction, beck2017first, doi:10.1137/18M121160X}. Let $V$ be a real Hilbert space equipped with the inner product $\langle \cdot, \cdot \rangle$ and induced norm $\| \cdot \|$.

For an extended real-valued function $f: V \to (-\infty, \infty]$, the \textit{domain} of $f$ is the set $\operatorname{dom} (f) = \{ x \in V \mid f(x) < \infty \}$. The \textit{epigraph} of $f$ is defined by $\operatorname{epi} (f) = \{ (x, y) \mid f(x) \leq y, x \in V, y \in \mathbb{R} \}$. A function $f$ is called \textit{proper} if $\operatorname{dom} (f) \neq \emptyset$. A function $f$ is \textit{closed} if $\operatorname{epi} (f)$ is closed. 

\begin{definition}[convex functions~\cite{beck2017first}]
$f: V \to (-\infty, \infty]$ is convex if $\operatorname{dom} (f)$ is convex, and for any $x, y \in \operatorname{dom} (f)$ and $\theta \in [0, 1]$,
    \begin{equation} \label{eq:convex}
        f(\theta x + (1 - \theta) y) \leq \theta f(x) + (1 - \theta) f(y).
    \end{equation}
\end{definition}

\begin{definition}[$\mu$-weakly convex functions~\cite{doi:10.1137/18M121160X}] \label{def:mu-weakly-convex}
    A function $f: V \to (-\infty, \infty]$ is $\mu$-weakly convex, $\mu \ge 0$, if the function $x \mapsto f(x) + \frac{\mu}{2} \| x \|^2$ is convex.
\end{definition}

The \textit{subdifferential} of a proper function $f$ at $x \in V$ is the set-valued mapping
\begin{equation*}
    \partial f(x):=\{u\in V \mid \forall y\in V,\ \langle y-x,u\rangle+f(x)\le f(y)\}.
\end{equation*}

For a set-valued operator $\operatorname{D}: V \to 2^V$, its \textit{graph} is $\operatorname{gra} \operatorname{D} = \{ (x, u) \in V \times V \mid u \in \operatorname{D}x\}$. The \textit{resolvent} of $\operatorname{D}$ with $\tau > 0$ is defined as $\operatorname{J}_{\tau \operatorname{D}} := \left(\operatorname{Id} + \tau \operatorname{D}\right)^{-1}$, where $\operatorname{Id}$ is the identity operator. In particular, an operator $\operatorname{D}: V \to 2^V$ such that, for every $x \in V$, $\operatorname{D}x$ is a singleton, then $\operatorname{D}$ is said to be \textit{(at most) single-valued}. In this paper, we restrict attention to the single-valued operators. 

\begin{definition}[firmly nonexpansive~\cite{bauschke2017correction}] \label{def:firm}
    An operator $\operatorname{D}: V \to V$ is \textit{firmly nonexpansive} if, for any $x, y \in V$, $\| \operatorname{D}x - \operatorname{D}y \|^2 \leq \langle x - y, \operatorname{D}x - \operatorname{D}y \rangle$. 
\end{definition}

\begin{definition}[nonexpansive~\cite{bauschke2017correction}] \label{def:ne}
    An operator $\operatorname{D}: V \to V$ is \textit{nonexpansive} if, for any $x, y \in V$, $\| \operatorname{D} x - \operatorname{D} y \| \leq \| x - y \|$.
\end{definition}

\begin{definition}[pseudo-contractive (PC)~\cite{bauschke2017correction}]\label{def:PC}
    An operator $\operatorname{D}:V \to V$ is \textit{pseudo-contractive} with parameter $k$, $k \in [0,1]$, if for any $x,y\in V$,
    \begin{equation}\label{eq:k-SPC}
        \| \operatorname{D}x-\operatorname{D}y\|^2 \le \|x-y\|^2 + k \|(\operatorname{Id}-\operatorname{D})x-(\operatorname{Id}-\operatorname{D})y\|^2.
    \end{equation}
    When $k \in (0, 1)$, $\operatorname{D}$ is $k$-strictly pseudo-contractive ($k$-SPC). 
\end{definition}

Note that nonexpansiveness is a special case of pseudo-contractivity with $k = 0$. The relationships among these operator properties are summarized as follows~\cite{AHQS}:
\begin{equation} \label{eq:relation_operators}
    \text{firmly nonexpansive} \;\implies\; \text{nonexpansive} \;\implies\; \text{PC},
\end{equation}
indicating that pseudo-contractivity imposes the weakest restriction among them.

\begin{definition}[$\beta$-cocoercive~\cite{bauschke2017correction}] \label{def:beta-coco} 
    An operator $\operatorname{D}:V \to V$ is \textit{$\beta$-cocoercive} for $\beta \ge 0$, if for any $x,y \in V$, 
    \begin{equation} \label{eq:beta-coco}
        \langle \operatorname{D} x-\operatorname{D} y,x - y\rangle \ge \beta \|\operatorname{D}x-\operatorname{D}y\|^2.
    \end{equation}
    If $\beta = 0$, the operator $\operatorname{D}$ is said to be \textit{monotone}. 
\end{definition}

Specifically, an operator $\operatorname{D}$ is \textit{maximally monotone} if there exists no monotone operator $\operatorname{A}$ such that $\operatorname{gra} \operatorname{A}$ properly contains $\operatorname{gra} \operatorname{D}$. 
\begin{lemma}[Theorem~20.25 and Example~23.3 in~\cite{bauschke2017correction}]
\label{lem:CCP-subdifferential}
    Let $f: V \to (-\infty, +\infty]$ be a proper, closed, and convex function.  
    Then its subdifferential operator $\partial f$ is maximally monotone. Moreover, for any $\tau > 0$, the resolvent of $\tau \partial f$ coincides with the proximal operator of $\tau f$, that is, for any $x \in V$,
    \begin{equation}
    \begin{aligned}
        \operatorname{J}_{\tau \partial f}(x)
        =& \operatorname{Prox}_{\tau f}(x) \\
        :=& \arg\min_{z \in V} \left\{ f(z) + \frac{1}{2\tau} \| z - x \|^2 \right\}.
    \end{aligned}
    \end{equation}
\end{lemma}

\subsection{Monotone inclusion problem}
A monotone inclusion problem (MIP)~\cite{Lee2025DYSSRG} is generally defined as finding $x \in V$, such that $0\in \operatorname{A}x$, where $\operatorname{A}$ is a maximal monotone operator. This framework encompasses classical (convex) minimization, variational inequalities, and saddle point problems, offering a unified and robust approach to convergence analysis that surpasses standard optimization methods in both flexibility and theoretical rigor.  

To the best of our knowledge, recent multi-prior tensor-recovery methods~\cite{zhao2020deep, zhao2021tensor, liu2025tensor, liang2024fixed, 10756233, HU2025112612} remain rooted in the classical optimization paradigm and typically treat the denoiser $\operatorname{D}_\sigma$ as a proximal map, solving the resulting problems via ADMM-type schemes. As discussed in section~\ref{subsec:related_works}, their convergence guarantees are largely empirical or rely on restrictive and often impractical assumptions about the denoiser~\cite{hurault2022proximal, wei2025learning}. Liang et al.~\cite{liang2024fixed} established fixed-point convergence for a multi-block PnP-ADMM under a bounded-denoiser assumption, while Hu et al.~\cite{HU2025112612} proved fixed-point convergence for PnP-ADMM under a nonexpansive assumption on deep neural networks. However, the resulting fixed point generally does not correspond to the solution of any optimization objective, which limits interpretability. Fig.~\ref{fig:VennMIP} illustrates the distinction between existing multi-prior methods and our proposed approach: by formulating the problem within the monotone inclusion framework, we treat priors as general operators rather than proximal mappings, thereby imposing substantially milder constraint on $\operatorname{D}_\sigma$. In this paper, we focus on the following MIP with three operators:
\begin{equation} \label{eq:MIP}
    \text{find  } x \in V, \text{ such that } 0\in \operatorname{A} x + \operatorname{B} x + \operatorname{C} x,
\end{equation}
where $\operatorname{A},\operatorname{B},\operatorname{C}$ are maximally monotone operators, and $\operatorname{C}$ is additionally $\beta$-cocoercive. 

\section{Method} \label{sec:tc_mip}
In this work, we utilize the MIP framework to address the tensor completion problem with multi-priors, constraining Eq.~\eqref{eq:MIP} to the tensor space $\mathbb{V}$:
\begin{equation} \label{eq:tensor-MIP}
    \text{find  } \mathcal{X} \in \mathbb{V} \text{ such that } 0 \in \operatorname{A} \mathcal{X} + \operatorname{B} \mathcal{X} + \operatorname{C} \mathcal{X}.
\end{equation}
According to Lemma~\ref{lem:CCP-subdifferential}, the operators $\operatorname{A}$ and $\operatorname{B}$ can be defined as the subdifferentials of convex functions on $\mathbb{V}$. Typically, $\operatorname{A}$ corresponds to the subdifferential of a convex data fidelity term. In the tensor completion setting, we define $\operatorname{A} = \partial \delta_{\mathcal{Y}, \Omega}$, where $\mathcal{Y} \in \mathbb{V}$ denotes the observed tensor with index set $\Omega$, and $\delta_{\mathcal{Y}, \Omega}$ is the indicator function enforcing data consistency:
\begin{equation*} \label{eq:delta}
    \delta_{\mathcal{Y}, \Omega} \left( \mathcal{X} \right) = \left\{
    \begin{array}{ll}
        0, & \mathcal{P}_{\Omega} \left( \mathcal{X} \right) = \mathcal{P}_{\Omega} \left( \mathcal{Y} \right), \\
        +\infty, & \text{otherwise.}
    \end{array}\right.
\end{equation*}
Within the formulation of Eq.~\eqref{eq:tensor-MIP}, any convex low-rank prior can be seamlessly integrated by assigning its subdifferential to $\operatorname{B}$. To further improve completion quality under highly undersampled observations, we extend the state-of-the-art tensor correlated total variation (t-CTV) prior to a more general and flexible form, termed the \textit{generalized t-CTV} (GTCTV), and set $\operatorname{B}$ accordingly in our model. To better preserve fine local details, we additionally incorporate \textit{deep pseudo-contractive} (DPC) denoisers as $\operatorname{C}$. The following sections present the formulations of GTCTV and DPC, followed by the proposed monotone-inclusion-based tensor completion model and its corresponding algorithm, along with a rigorous global convergence analysis.

\subsection{Generalized tensor correlated total variation} \label{subsec:GTCTV}
To jointly promote the low-rank and smooth features in tensors, Wang et al.~\cite{10078018} introduced t-CTV, demonstrating superior recovery performance. Building on this idea, we extend TNN~\cite{9730793} with weakly convex penalty to formulate a more flexible GTCTV, thereby further enhancing low‑rank characteristics in the gradient domain.

\begin{definition}[tensor $f$-penalty] \label{def:tensor-f-penalty}
    Let $f: \mathbb{R}_{\geq 0} \to \mathbb{R}_{\geq 0}$ be a $\mu$-weakly convex penalty such that $\operatorname{Prox}_{\eta f}$ is non-decreasing for any $\eta > 0$. Given an invertible linear transform~$\mathfrak{L}$, the tensor $f$-penalty of $\mathcal{A}$ is defined as:
    \begin{equation*}
        \begin{aligned}
            \| \mathcal{A} \|_{f, \mathfrak{L}} &= \frac{1}{l^2} \sum_{i_3 \dots i_N} \left\| \mathcal{A}^{(i_3 \dots i_N)}_{\mathfrak{L}} \right\|_{f} \\
            &= \frac{1}{l^2} \sum_{i_3 \dots i_N} \sum_{j=1}^r f \left( \sigma_j \left( \mathcal{A}^{(i_3 \dots i_N)}_{\mathfrak{L}} \right) \right), \\
        \end{aligned}
    \end{equation*}
    where $r=\min \{n_1, n_2\}$, and $l$ is the composite scale factor corresponded to $\mathfrak{L}$.
\end{definition}

Leveraging the flexibility of the tensor $f$-penalty, we can employ a wide class of weakly convex regularizers, such as MCP~\cite{zhangnearly} and SCAD~\cite{fanvariable}, whose proximal solutions are non-decreasing, to more aggressively promote low-rankness in each gradient tensor. The resulting GTCTV is then obtained by averaging these $f$-penalties over a predefined set of directional gradients. 
\begin{definition}[Generalized Tensor Correlated Total Variation (GTCTV)] \label{def:GTCTV}
    Let $\mathcal{A} \in \mathbb{V}$, $\Gamma \subseteq \{ 1, 2, \dots, N\}$ be a predefined set of gradient directions, and let $\gamma = \sharp \{ \Gamma \}$ denote its cardinality, with $\nabla_d$ be the linear gradient operator along the $d$-th mode defined in Table~\ref{table:notations}. Then the GTCTV of $\mathcal{A}$ with respect to a regularizer $f$ and transform $\mathfrak{L}$ is 
    \begin{equation} \label{eq:GTCTV}
        \| \mathcal{A} \|_{\operatorname{GTCTV}} := \frac{1}{\gamma} \sum_{d \in \Gamma} \left\| \nabla_d \mathcal{A} \right\|_{f, \mathfrak{L}}.
    \end{equation}
\end{definition}

We remark that when $f(x) = |x|$, the proposed GTCTV prior reduces to the original t-CTV. To facilitate its incorporation into the MIP formulation in Eq.~\eqref{eq:tensor-MIP}, we further provide a rigorous analysis of the weak convexity of GTCTV for a general $\mu$-weakly convex function $f$, as established in Lemma~\ref{lem:mu0-gtctv}.
\begin{lemma}[Proof in section~S.\Rmnum{1}.B of the supplement] \label{lem:mu0-gtctv}
    Let $f$ be a $\mu$-weakly convex function. Then, the prior $\|\cdot\|_{\text{GTCTV}}$ is $4 \mu$-weakly convex on $\mathbb{V}$.
\end{lemma}
    
\subsection{Deep pseudo-contractive denoisers} \label{subsec:dpc}
Deep PnP denoisers have proven effective for tensor restoration, notably in MDI inpainting~\cite{GUO2023171} and traffic data completion~\cite{10756233}. These tasks employ denoisers to solve proximal subproblems within PnP methods~\cite{zhang2021plug}. However, many PnP methods rely on strong theoretical assumptions that are difficult to satisfy for deep denoisers~\cite{hurault2022proximal, wei2025learning}. 

To address this, Wei et al.~\cite{AHQS} introduced DPC denoisers, employing a loss function to enforce approximate pseudo-contractive properties. This approach, based on less stringent assumptions, improves PnP-based image restoration quality without imposing architectural restrictions on the network. The weaker property in Eq.~\eqref{eq:relation_operators} imply fewer constraint on the operator. As noted in~\cite{AHQS}, when the denoiser $\operatorname{D}_{\sigma}$ satisfies weaker assumption, it empirically exhibits improved denoising performance, enhancing the efficacy of PnP iterative frameworks for image restoration.

Specifically, Wei et al.~\cite{AHQS} developed a loss function that encourages the denoiser $\operatorname{D}_{\sigma}$ to be $k$-SPC ($k \in (0, 1)$) for noisy images $\mathbf{X} \in \mathbb{R}^{n_1 \times n_2 \times n_3}$, where $n_3 = 1$ or $3$. In our application, which focuses on order-$N$ tensors $\mathcal{X} \in \mathbb{V}$, we apply $\operatorname{D}_{\sigma}$ considering the following two cases:
\begin{itemize}
    \item \textbf{Case 1:} If $n_3 = 1$ or $3$, we define $\operatorname{D}_{\sigma} (\mathcal{X})$ slice-wise as 
    \begin{equation*}
        \left[\operatorname{D}_{\sigma} (\mathcal{X})\right]^{(:, i_4\dots i_N)} = \operatorname{D}_{\sigma} (\mathcal{X}^{(:, i_4\dots i_N)}),
    \end{equation*}
    where $\mathcal{X}^{(:, i_4\dots i_N)} := \mathcal{X}(:, :, :, i_4, \dots, i_N)$.
    \item \textbf{Case 2:} If $n_3 \neq 1$ and $3$, we expand the dimensions to $N+1$ and set $n_3 = 1$. For example, a $256 \times 256 \times 31$ MSI is reshaped to $256 \times 256 \times 1 \times 31$, thereby reducing it to \textbf{Case 1} without affecting the Frobenius norm.
\end{itemize}

Therefore, for any $\mathcal{X}$, $\mathcal{Y} \in \mathbb{V}$, 
\begin{equation*}
    \begin{aligned}
        & \| \operatorname{D}_{\sigma} \left( \mathcal{X} \right) - \operatorname{D}_{\sigma} \left( \mathcal{Y} \right) \|_F^2 \\
        =& \sum_{i_4, \dots, i_N} \left\| \operatorname{D}_{\sigma} \left( \mathcal{X}^{(:, i_4\dots i_N)} \right) - \operatorname{D}_{\sigma} \left( \mathcal{Y}^{(:, i_4\dots i_N)} \right) \right\|_F^2 \\
        \leq & \sum_{i_4, \dots, i_N} \left( \left\| \mathcal{X}^{(:, i_4\dots i_N)} - \mathcal{Y}^{(:, i_4\dots i_N)} \right\|_F^2 \right.\\
        &\left. + k \left\| \left(\operatorname{Id} - \operatorname{D}_{\sigma} \right) \left(\mathcal{X}^{(:, i_4\dots i_N)}\right) - \left(\operatorname{Id} - \operatorname{D}_{\sigma} \right) \left(\mathcal{Y}^{(:, i_4\dots i_N)}\right) \right\|_F^2 \right) \\
        =& \left\| \mathcal{X} - \mathcal{Y} \right\|_F^2 + k \left\| \left(\operatorname{Id} - \operatorname{D}_{\sigma} \right) \left(\mathcal{X}\right) - \left(\operatorname{Id} - \operatorname{D}_{\sigma} \right) \left(\mathcal{Y}\right) \right\|_F^2, \\
    \end{aligned}
\end{equation*}
Thus, $\operatorname{D}_{\sigma}$ is $k$-SPC on $\mathbb{V}$ by Definition~\ref{def:PC}. The link between pseudo-contractivity and $\beta$-cocoercivity is formalized in Lemma~\ref{lem:spc and coco}. This result is crucial for formulating our tensor completion model within the MIP framework in Eq.~\eqref{eq:MIP}.
\begin{lemma}[Proof in section~S.\Rmnum{1}.C of the supplement] \label{lem:spc and coco}
    $\operatorname{D}:V \to V$ be a PC operator with $k \in [0,1]$, if and only if, $\operatorname{Id} - \operatorname{D}$ is $\beta$-cocoercive with $\beta = \frac{1-k}{2}$. 
\end{lemma}

\begin{remark} \label{remark:scale-cocoercive}
    For a $\beta_0$-cocoercive operator $\operatorname{A}$, $\alpha > 0$ and any $x, y \in V$,
    \begin{equation*}
        \begin{aligned}
            &\langle \alpha \operatorname{A} x - \alpha \operatorname{A} y, x - y \rangle = \alpha \langle \operatorname{A} x - \operatorname{A} y, x - y \rangle \\
            \geq &\alpha \beta_0 \| \operatorname{A} x - \operatorname{A} y \|^2 = \frac{\beta_0}{\alpha} \| \alpha \operatorname{A} x - \alpha \operatorname{A} y \|^2.
        \end{aligned}
    \end{equation*}
    Thus, $\alpha \operatorname{A}$ is $\frac{\beta_0}{\alpha}$-cocoercive.
\end{remark}

\subsection{Tensor completion via MIP}
Drawing upon the aforementioned priors, we address tensor completion within the MIP framework of Eq.~\eqref{eq:tensor-MIP} by setting
\begin{equation} \label{eq:concrete_set}
    \operatorname{A} = \partial \delta_{\mathcal{Y}, \Omega}, \, \operatorname{B} = \partial \left(\left\| \cdot \right\|_{\text{GTCTV}} + 2 \mu \left\|\cdot\right\|_F^2 \right), \, \operatorname{C} = \alpha \left( \operatorname{Id} - \operatorname{D}_{\sigma}\right).
\end{equation}
Here, $\mu$ is the weak convexity parameter of the base function $f$, and the GTCTV prior is $4 \mu$-weakly convex, as shown in Lemma~\ref{lem:mu0-gtctv}. Additionally, $\alpha > 0$ is a weighting factor, and $\operatorname{D}_\sigma$ is a deep $k$-SPC Gaussian denoiser with denoising strength $\sigma$. This yields the following concrete tensor completion model: find $\mathcal{X} \in \mathbb{V}$, such that
\begin{equation} \label{eq:Proposed-MIP}
    \begin{aligned}
        \mathcal{O} \in & \partial \delta_{\mathcal{Y}, \Omega} \left( \mathcal{X} \right) + \partial \left(\left\| \mathcal{X} \right\|_{\text{GTCTV}} + 2 \mu \left\|\mathcal{X}\right\|_F^2 \right) \\
        & + \alpha \left( \operatorname{Id} - \operatorname{D}_{\sigma}\right) \left( \mathcal{X} \right). 
    \end{aligned}
\end{equation}

Consequently, $\operatorname{Id} - \operatorname{D}_{\sigma}$ outputs the predicted noise. Distinct from the traditional PnP paradigm, which incorporates $\operatorname{D}_\sigma$ in a backward fashion (e.g., as a proximal operator or resolvent), our methodology utilizes $\operatorname{D}_\sigma$ in a forward manner. This approach aligns with frameworks such as RED~\cite{romano2017little}, and diffusion-based techniques~\cite{chung2023diffusion}. 
 
Moreover, assuming $R_{\text{global}}\left(\mathcal{X}\right) = \left\| \mathcal{X} \right\|_{\text{GTCTV}} + 2 \mu \left\|\mathcal{X}\right\|_F^2$ and $R_{\text{local}}$ is an implicit function whose gradient is $\operatorname{Id} - \operatorname{D}_{\sigma}$, the MIP in Eq.~\eqref{eq:Proposed-MIP} can be interpreted as the first-order optimality condition of the classical regularized optimization for tensor completion:
\begin{equation} \label{eq:OP}
    \min_{\mathcal{X} \in \mathbb{V}} \delta_{\mathcal{Y}, \Omega} \left(\mathcal{X}\right) + R_{\text{global}}\left(\mathcal{X}\right) + R_{\text{local}}\left(\mathcal{X}\right).
\end{equation}

However, to the best of our knowledge, designing iterative schemes that solve Eq.~\eqref{eq:OP} with guaranteed global convergence remains a significant challenge, highlighting the robustness and appeal of the monotone inclusion framework. Moreover, Wei et al.~\cite{wei2025learning} show that a well-defined $R_{\text{local}}$ exists only if the deep denoiser is conservative. Imposing such a requirement introduces additional constraints on the denoiser, which may limit its effectiveness and ultimately degrade overall recovery performance.

\subsection{The proposed algorithm: GTCTV-DPC}
To solve the general MIP in Eq.~\eqref{eq:MIP}, Davis and Yin~\cite{davis2017three} proposed the well-known Davis-Yin splitting (DYS) method. They reframe the MIP as a fixed-point problem:
\begin{equation*}\label{eq:FPP}
    \text{find  } z \in V \text{, such that } z=\operatorname{T}z,
\end{equation*}
where $\operatorname{T}$ with stepsize $\tau$ is defined as
\begin{equation}\label{eq:T in DYS}
    \operatorname{T}:=\operatorname{Id}-\operatorname{J}_{\tau \operatorname{B}}+\operatorname{J}_{\tau \operatorname{A}}\circ( 2\operatorname{J}_{\tau \operatorname{B}}-\operatorname{Id}-\tau \operatorname{C} \circ \operatorname{J}_{\tau \operatorname{B}}).
\end{equation}
The detailed DYS algorithm, which adopts a Krasnoselskii-Mann (KM)–type iteration~\cite{hicks1977mann, rafiq2007mann}, is summarized in Algorithm~\ref{alg:DYS}.
\begin{algorithm}
\caption{DYS Algorithm for Solving the MIP in Eq.~\eqref{eq:MIP}.} \label{alg:DYS}
\begin{algorithmic}
    \State \textbf{Input: } $z_0$, $\tau$, $N_{\text{max}}$, $\{\lambda_t\}_{t=0}^{N_{\text{max}}}$;
    \For{$t=0:N_{\text{max}}-1$}
        \State $x_{t+1}^{\operatorname{B}} = \operatorname{J}_{\tau \operatorname{B}} \left( z_t \right) $;
        \State $x_{t+1}^{\operatorname{A}} = \operatorname{J}_{\tau \operatorname{A}} \left( 2 x_{t+1}^{\operatorname{B}} - z_t - \tau \operatorname{C} \left( x_{t+1}^{\operatorname{B}} \right) \right)$; 
        \State \AlgComment{$z_{t+1} = (1 - \lambda_t) z_t + \lambda_t \operatorname{T} z_t$.}
        \State $z_{t+1} = z_t + \lambda_t \left( x_{t+1}^{\operatorname{A}} - x_{t+1}^{\operatorname{B}}\right)$. 
    \EndFor
    \State Return $x_{t+1}^{\operatorname{A}}$.
\end{algorithmic}
\end{algorithm}

We employ the DYS method in Algorithm~\ref{alg:DYS} to solve the proposed MIP in Eq.~\eqref{eq:Proposed-MIP}. We begin by outlining the computation of the relevant resolvent operators with stepsize $\tau$, followed by the concrete algorithm. 

\subsubsection{Resolvent of $\tau \partial \delta_{\mathcal{Y}, \Omega}$}
Since $\delta_{\mathcal{Y}, \Omega}$ is convex, $\partial \delta_{\mathcal{Y}, \Omega}$ is maximally monotone by Lemma~\ref{lem:CCP-subdifferential}. Consequently, the resolvent of $\tau \partial \delta_{\mathcal{Y}, \Omega}$ is the proximal operator of $\delta_{\mathcal{Y}, \Omega}$: 
\begin{equation} \label{eq:resolvent-delta}
     \operatorname{J}_{\tau \partial \delta_{\mathcal{Y}, \Omega}} \left( \mathcal{X} \right)  = \operatorname{Prox}_{\tau \delta_{\mathcal{Y}, \Omega}} \left( \mathcal{X} \right) 
     = \mathcal{P}_{\Omega} \left(\mathcal{Y}\right) + \mathcal{P}_{\Omega^{\perp}} \left(\mathcal{X}\right),
\end{equation}
where $\Omega^{\perp}$ denotes the complement of $\Omega$. 

\subsubsection{Resolvent of $\tau \partial ( \| \cdot\|_{\text{GTCTV}} + 2 \mu \|\cdot\|_F^2 )$}
Given the convexity of $\| \cdot\|_{\text{GTCTV}} + 2 \mu \|\cdot\|_F^2$ as proven in Lemma~\ref{lem:mu0-gtctv}, the resolvent operator reduces to the proximal operator associated with $\tau \left(\| \cdot\|_{\text{GTCTV}} + 2 \mu \|\cdot\|_F^2 \right)$: 
\begin{equation*} \label{eq:resolvent-gtctv}
    \begin{aligned}
        & \operatorname{J}_{\tau \partial \left( \| \cdot\|_{\text{GTCTV}} + 2 \mu \|\cdot\|_F^2 \right)} \left( \mathcal{X} \right) = \operatorname{Prox}_{\tau \left(\| \cdot\|_{\text{GTCTV}} + 2 \mu \|\cdot\|_F^2 \right)} \left( \mathcal{X} \right) \\
        = & \arg\min\limits_{\mathcal{M} \in \mathbb{V}} \left( \left\| \mathcal{M} \right\|_{\text{GTCTV}} + 2 \mu \left\|\mathcal{M}\right\|_F^2  + \frac{1}{2 \tau } \left\| \mathcal{M} - \mathcal{X} \right\|_F^2 \right). \\ 
    \end{aligned}
\end{equation*}

To separate the difference operation $\nabla_d (\cdot)$, we introduce auxiliary variables $\mathcal{G}_d$ and employ the ADMM~\cite{boyd2011distributed} to solve the reformulated subproblem:
\begin{equation*} \label{eq:admm-gtctv}
    \begin{aligned}
        \min_{\substack{\mathcal{M} \in \mathbb{V}, \\\mathcal{G}_d, d \in \Gamma}} &\frac{1}{\gamma} \sum_{d \in \Gamma} \left\| \mathcal{G}_d \right\|_{f, \mathfrak{L}} + 2 \mu \left\|\mathcal{M}\right\|_F^2 + \frac{1}{2\tau} \left\| \mathcal{M} - \mathcal{X} \right\|_F^2 \\
        \text{ s.t. } &\mathcal{G}_d = \nabla_d\mathcal{M}, d \in \Gamma.
    \end{aligned}
\end{equation*}
The augmented Lagrangian function is
\begin{equation} \label{eq:lagrange-gtctv}
    \begin{aligned}
        & \mathcal{L} \left( \mathcal{M}, \left\{\mathcal{G}_d, d\in \Gamma\right\},  \left\{\mathcal{B}_d, d\in \Gamma\right\} \right) \\
        = & \sum_{d \in \Gamma} \left( \frac{1}{\gamma} \left\| \mathcal{G}_d \right\|_{f, \mathfrak{L}}  + \frac{\rho_t}{2} \left\| \nabla_d \mathcal{M} - \mathcal{G}_d + \frac{\mathcal{B}_d}{\rho_t}\right\|_F^2 \right) \\
        &+ 2 \mu \left\|\mathcal{M}\right\|_F^2 + \frac{1}{2\tau} \left\| \mathcal{M} - \mathcal{X} \right\|_F^2, \\
    \end{aligned}
\end{equation}
where $\rho_t > 0$ is a penalty parameter, and $\mathcal{B}_d$ is the Lagrange multiplier. We describe how to solve the subproblems for each variable as follows:

$\bullet$ \textit{Updating $\mathcal{M}^{t+1}$}: Following~\cite{doi:10.1137/080724265, 10078018}, we compute the derivative of Eq.~\eqref{eq:lagrange-gtctv} with respect to $\mathcal{M}$:
    \begin{equation} \label{eq:linear_system_m}
        \begin{aligned}
            &\left( \tau \rho_t \sum_{d \in \Gamma} \nabla_d^{\top} \nabla_d + (4 \tau \mu + 1) \operatorname{Id} \right) \mathcal{M} \\
            =& \tau \sum_{d \in \Gamma}\nabla_d^{\top} (\rho_t \mathcal{G}_d - \mathcal{B}_d) + \mathcal{X}. 
        \end{aligned}
    \end{equation}
    
    Applying multi-dimensional FFT to solve Eq.~\eqref{eq:linear_system_m} yields the optimal solution for $\mathcal{M}^{t+1}$: 
    \begin{equation} \label{eq:M-solution-gtctv}
        \mathcal{M}^{t+1} = \operatorname{F}^{-1} \left( \frac{\operatorname{F} \left( \mathcal{X}\right) + \tau  \sum_{d \in \Gamma}\operatorname{F} \left( \mathcal{D}_d\right)^* \odot \operatorname{F} \left( \rho_t \mathcal{G}_d^{t} - \mathcal{B}_d^{t}\right)}{(4 \tau \mu + 1)\mathbf{1} + \tau \rho_t \sum_{d\in\Gamma} \operatorname{F} \left( \mathcal{D}_d\right)^* \odot \operatorname{F}\left( \mathcal{D}_d\right)}\right),
    \end{equation}
    where $\mathcal{D}_d$ denotes the difference tensor of $\nabla_d$, $\operatorname{F}$ represents the multi-dimensional FFT along all modes, $\mathbf{1}$ is a tensor with all entries equal to $1$, $\odot$ indicates componentwise multiplication, and the division is performed componentwise as well.

$\bullet$ \textit{Updating $\mathcal{G}_d^{t+1}$}: For each $d\in\Gamma$, isolating the terms involving $\mathcal{G}_d$ in Eq.~\eqref{eq:lagrange-gtctv} results in the subproblem:
    \begin{equation} \label{eq:G-subproblem-gtctv}
        \mathcal{G}_d^{t+1} = \arg\min\limits_{\mathcal{G}_d \in \mathbb{V}} \frac{1}{\gamma} \left\| \mathcal{G}_d \right\|_{f, \mathfrak{L}}  + \frac{\rho_t}{2} \left\| \mathcal{G}_d - \left( \nabla_d \mathcal{M}^{t+1} + \frac{\mathcal{B}_d^t}{\rho_t} \right) \right\|_F^2. 
    \end{equation}

    To solve the subproblem in Eq.~\eqref{eq:G-subproblem-gtctv}, we provide the proximal solution for the tensor $f$-penalty in Lemma~\ref{lem:tensor-f-shrinkage}.
    \begin{lemma}[Proof in section~S.\Rmnum{1}.D of the supplement] \label{lem:tensor-f-shrinkage}
    Given $\| \cdot \|_{f, \mathfrak{L}}$ as defined in Definition~\ref{def:tensor-f-penalty}, and a tensor $\mathcal{T} \in \mathbb{V}$ with t-SVD $\mathcal{T} = \mathcal{U} *_{\mathfrak{L}} \mathcal{S} *_{\mathfrak{L}} \mathcal{V}^{*}$, the solution to the proximal problem 
    \begin{equation*}
        \mathcal{G}_* = \arg \min_{\mathcal{G} \in \mathbb{V}} \| \mathcal{G} \|_{f, \mathfrak{L}} + \frac{1}{2\eta} \| \mathcal{G} - \mathcal{T}\|_F^2
    \end{equation*}
    is given by $\mathcal{G}_* = \operatorname{t-SVF}_{\eta f} \left(\mathcal{T} \right) := \mathcal{U} *_{\mathfrak{L}} \mathcal{S}_{\eta f} *_{\mathfrak{L}} \mathcal{V}^{*}$, where $[\mathfrak{L} (\mathcal{S}_{\eta f})]^{(i_3 \dots i_N)} = \operatorname{diag} (\boldsymbol{\sigma})$, and $\boldsymbol{\sigma}_i = \operatorname{Prox}_{\eta f} \left( \left[ \mathcal{S}_{\mathfrak{L}} \right]_{i, i, i_3, \dots, i_N} \right)$ for $i = 1, 2, \dots, r$. 
    \end{lemma}

    Applying Lemma~\ref{lem:tensor-f-shrinkage} to Eq.~\eqref{eq:G-subproblem-gtctv}, we obtain: 
    \begin{equation} \label{eq:G-solution-gtctv}
        \mathcal{G}_d^{t+1} = \operatorname{t-SVF}_{\frac{1}{\gamma \rho_t} f} \left( \nabla_d \mathcal{M}^{t+1} + \frac{\mathcal{B}_d^t}{\rho_t} \right).
    \end{equation}

\begin{figure*}[htbp!]
    \centering
    \includegraphics[width=.92\linewidth]{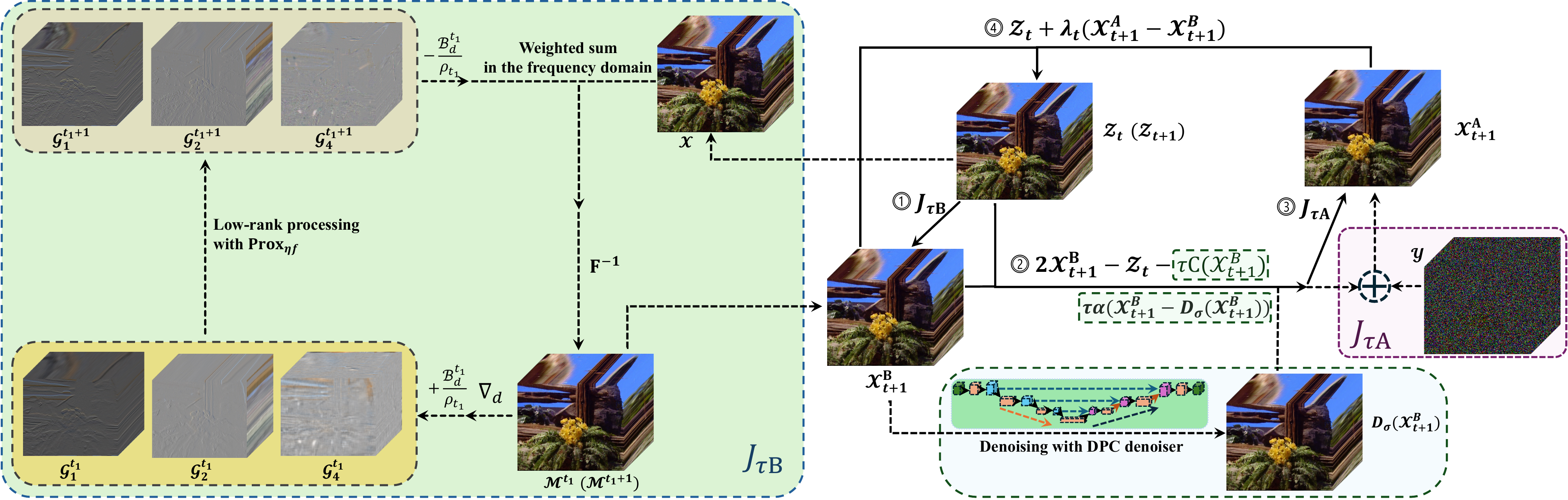}
    \caption{Iterative flowchart of GTCTV-DPC, illustrated using a color video example. In the flowchart, the symbol $\bigoplus$ represents the operation in Eq.~\eqref{eq:resolvent-delta}, and $\operatorname{F}^{-1}$ denotes the inverse multi-dimensional FFT in Eq.~\eqref{eq:M-solution-gtctv}.}
    \label{fig:algorithm_review}
\end{figure*}

\subsubsection{GTCTV-DPC}
The entire algorithm is summarized in Algorithm~\ref{alg:main}, and Fig.~\ref{fig:algorithm_review} illustrates the main iterative flow of GTCTV-DPC.
\begin{algorithm}[htbp!]
\caption{GTCTV-DPC for Tensor Completion via MIP.} \label{alg:main}
\begin{algorithmic}[1]
    \State \textbf{Input: } Observation $\mathcal{Y}$, the step size $\tau$, $f$, $\alpha$, $\rho_0$, $k$-SPC deep denoiser $\operatorname{D}_{\sigma_t}$, $\sigma_0$, $\nu$, $\varepsilon$, $N_{\text{in}}$, $N_{\text{max}}$, and $\{ \lambda_t\}_{t=0}^{N_{\text{max}}}$.
    \State \textbf{Initialize: } $\mathcal{Z}_{0}=\mathcal{Y}$, $\mathcal{G}_d^{t_1+1} = \nabla_d \mathcal{Z}_{0}$, $\mathcal{B}_d^{t_1 + 1} = \mathcal{O}$.
    \For{$t=0:N_{\text{max}}-1$}
        \State \AlgComment{$\mathcal{X}_{t+1}^{\operatorname{B}} = \operatorname{J}_{\tau \operatorname{B}} \left( \mathcal{Z}_t \right) $;}
        \State Let $\mathcal{X} = \mathcal{Z}_t$, $\mathcal{G}_d^{0} = \mathcal{G}_d^{t_1+1}$, $\mathcal{B}_d^{0} = \mathcal{B}_d^{t_1+1}$; 
        \For{$t_1 = 0 : N_{\text{in}} - 1$}
            \State Update $\mathcal{M}^{t_1 + 1}$ by Eq.~\eqref{eq:M-solution-gtctv}; 
            \State Update $\mathcal{G}_d^{t_1 + 1}$ by Eq.~\eqref{eq:G-solution-gtctv}, $d \in \Gamma$; 
            \State $\mathcal{B}_d^{t_1+1} = \mathcal{B}_d^{t_1} + \rho_{t_1} \left( \nabla_d \mathcal{M}^{t_1+1} - \mathcal{G}_d^{t_1+1}\right)$, $d \in \Gamma$;
            \State \AlgComment{For accelerating convergence.}
            \State Let $\rho_{t_1 + 1} = \min \{ \nu \rho_{t_1}, 10^{10} \}$; 
            \State $\epsilon_1 = \| \mathcal{M}^{t_1+1} - \mathcal{M}^{t_1} \|_F^2 \, / \, \| \mathcal{M}^{t_1} \|_F^2$; 
        \If{ $\epsilon_1 < \varepsilon$ }
            \State break.
        \EndIf 
        \EndFor
        \State $\mathcal{X}_{t+1}^{\operatorname{B}} = \mathcal{M}^{t_1 + 1}$; 
        \State \AlgComment{$\mathcal{X}_{t+1}^{\operatorname{A}} = \operatorname{J}_{\tau \operatorname{A}} \left( 2 \mathcal{X}_{t+1}^{\operatorname{B}} - z_t - \tau \operatorname{C} \left( \mathcal{X}_{t+1}^{\operatorname{B}} \right) \right)$;}
        \State $\mathcal{X}_{t+1}^{\operatorname{C}} = \mathcal{X}_{t+1}^{\operatorname{B}} - \operatorname{D}_{\sigma_t} \left( \mathcal{X}_{t+1}^{\operatorname{B}} \right)$; 
        \State $\mathcal{X}_{t+1}^{\operatorname{A}} = \mathcal{P}_{\Omega^{\bot}} \left( 2 \mathcal{X}_{t+1}^{\operatorname{B}} - \mathcal{Z}_{t} - \tau \alpha \mathcal{X}_{t+1}^{\operatorname{C}} \right) + \mathcal{P}_{\Omega} \left( \mathcal{Y} \right)$;
        \State \AlgComment{$\mathcal{Z}_{t+1}=(1-\lambda_t)\mathcal{Z}_{t}+\lambda_t \operatorname{T} \mathcal{Z}_{t}$;}
        \State $\mathcal{Z}_{t+1} = \mathcal{Z}_{t} + \lambda_t \left( \mathcal{X}_{t+1}^{\operatorname{A}} - \mathcal{X}_{t+1}^{\operatorname{B}} \right)$; 
        \State \AlgComment{For accelerating convergence.}
        \State Let $\sigma_{t+1} = \max \{ \sigma_t \, / \, \nu, 10^{-3} \}$; 
        \State $\epsilon = \| \mathcal{X}_{t+1}^{\operatorname{A}} - \mathcal{X}_{t}^{\operatorname{A}} \|_F^2 \, / \, \| \mathcal{X}_{t}^{\operatorname{A}} \|_F^2$; 
        \If{ $\epsilon < \varepsilon$ }
            \State break.
        \EndIf
    \EndFor
    \State Return $\mathcal{X}_{t+1}^{\operatorname{A}}$.
\end{algorithmic}
\end{algorithm}

\subsection{Convergence analysis} \label{sec:convergent_analysis}
In this section, we present a rigorous convergence analysis for GTCTV-DPC (Algorithm~\ref{alg:main}), formalized in Corollary~\ref{coro:convergence}. Starting from the general MIP on $V$ in Eq.~\eqref{eq:MIP}, we prove that the operator $\operatorname{T}$ in Eq.~\eqref{eq:T in DYS} is SPC for the extended range $\tau \in (0, 4 \beta)$, and we provide different admissible choices for the relaxation parameters $\{\lambda_t\}_{t \geq 0}$ in Theorem~\ref{the:convergence}, extending prior analyses~\cite{aragon2022direct, davis2017three}. 

\begin{theorem}[Proof in Section~S.\Rmnum{1}-E of the supplement] \label{the:convergence}
   Let $\operatorname{A}$ and $\operatorname{B}$ be maximally monotone, and $\operatorname{C}$ be $\beta$-cocoercive in Eq.~\eqref{eq:MIP}. Let $\operatorname{T}$ be the operator defined in Eq.~\eqref{eq:T in DYS} with stepsize $\tau\in (0, 4 \beta)$. Then, from any initial point $z_0 \in V$, the iteration $z_{t+1}=(1-\lambda_t)z_{t}+\lambda_t \operatorname{T} z_{t}$ converges weakly to a fixed point of $\operatorname{T}$, where $\{\lambda_t\}_{t\ge 0}$ is a real sequence in $[0, 1]$ satisfies 
   \begin{equation*}
       \text{(i) } \sum_{t=0}^\infty \lambda_t=\infty \text{ and (ii) }\sum_{t=0}^\infty \lambda_t^2 <\infty.
   \end{equation*}
    Moreover, $x_{t}^{\operatorname{B}}$ and $x_{t}^{\operatorname{A}}$ converge weakly to a solution of the MIP defined in Eq.~\eqref{eq:MIP}.
\end{theorem}

Based on Theorem~\ref{the:convergence}, we obtain the following Corollary~\ref{coro:convergence}, which guarantees the global convergence of Algorithm~\ref{alg:main}.

\begin{corollary}[Proof in Section~S.\Rmnum{1}-F of the supplement] \label{coro:convergence}
    Let $\operatorname{D}_\sigma$ be $k$-SPC with $k\in (0,1)$, $\tau\in (0,\frac{2-2k}{\alpha})$, and $\{\lambda_t\}_{t\ge 0}$ is a real sequence in $[0, 1]$ satisfy 
    \begin{equation*}
        \text{(i) } \sum_{t=0}^\infty \lambda_t=\infty \text{ and (ii) }\sum_{t=0}^\infty \lambda_t^2 <\infty.
    \end{equation*}
    Let $\operatorname{T}$ be the operator defined in Eq.~\eqref{eq:T in DYS}. Then, from any initial point $\mathcal{Z}_0 \in \mathbb{V}$, the iteration $\mathcal{Z}_{t+1}=(1-\lambda_t)\mathcal{Z}_{t}+\lambda_t \operatorname{T} \mathcal{Z}_{t}$ generated by Algorithm~\ref{alg:main} converges to a fixed point of $\operatorname{T}$. Moreover, $\mathcal{X}_{t}^{\operatorname{B}}$ and $\mathcal{X}_{t}^{\operatorname{A}}$ converge to a solution of the MIP defined in Eq.~\eqref{eq:Proposed-MIP}.
\end{corollary}

\section{Experiments} \label{sec:experiment}
In this section, we adopt a Bernoulli random sampling scheme across the entire tensor for all experiments, utilizing publicly available datasets. We evaluate the proposed GTCTV-DPC on two data types: multi-dimensional images and spatio-temporal traffic data. A comprehensive comparative study is conducted against state-of-the-art baselines. The source code is publicly available at the GitHub repository\footnote{https://github.com/peterchen96/TensorCompletionMIP.git}.

\subsection{Experimental settings} \label{subsec:expe_sets}
For MDI completion, we select twelve MSIs from the CAVE\footnote{https://cave.cs.columbia.edu/repository/Multispectral} dataset~\cite{5439932}, and eleven color videos from the YUV\footnote{http://trace.eas.asu.edu/yuv} dataset. Selected MSIs and color videos are shown in Figs~S.1 and S.2 from the supplement. To comprehensively evaluate GTCTV-DPC on MDI completion, we compare against the following methods: low-rank methods FTNN~\cite{9115254} and t-CTV~\cite{10078018}; deep learning methods HIR-Diff~\cite{pang2024hir}, LRTFR~\cite{10354352}, and DRO-TFF~\cite{Li_Zhang_Luo_Meng_2025}; multi-prior methods DP3LRTC~\cite{zhao2020deep}, and FBGND~\cite{liang2024fixed}. We evaluate MDI completion using Mean Peak Signal-to-Noise Ratio (MPSNR) and Mean Structural Similarity (MSSIM):
\begin{equation*}
\begin{aligned}
    \operatorname{MPSNR} &= \frac{1}{n_4} \sum_{i_4 = 1}^{n_4} \operatorname{PSNR} \left(\mathcal{X}_{\text{out}}^{(:, i_4)}, \mathcal{X}_{\text{ori}}^{(:, i_4)} \right), \\ \operatorname{MSSIM} &= \frac{1}{n_4} \sum_{i_4 = 1}^{n_4} \operatorname{SSIM} \left(\mathcal{X}_{\text{out}}^{(:, i_4)}, \mathcal{X}_{\text{ori}}^{(:, i_4)} \right), \\
\end{aligned}
\end{equation*}
where $\mathcal{X}^{(:, i_4)} := \mathcal{X}(:, :, :, i_4)$, and $\mathcal{X}_{\text{out}}$ and $\mathcal{X}_{\text{ori}}$ denote the completed and original tensors, respectively. For MSIs, we extend them to $256 \times 256 \times 1 \times 31$, treating them as the band number of grayscale images. Higher MPSNR and MSSIM values indicate better completion quality. 

For traffic data, we select three publicly available datasets from real-world transportation systems: Guangzhou\footnote{https://zenodo.org/records/1205229}, Seattle\footnote{https://github.com/zhiyongc/Seattle-Loop-Data}, and PeMS\footnote{https://people.eecs.berkeley.edu/\~varaiya/papers\_ps.dir/PeMSTutorial}. Each traffic dataset is structured either as a third-order tensor or as a time series matrix. For traffic data completion, we compare against: tensor/matrix factorization methods BATF~\cite{CHEN201966} and BTMF~\cite{9380704}; low-rank methods LRTC-TNN~\cite{CHEN2020102673} and LRTC-TSpN~\cite{NIE2022103737}; multi-prior methods LSTC-Tubal~\cite{CHEN2021103226} and LATC~\cite{chen2022low}; deep learning method LRTFR~\cite{10354352}. We evaluate traffic data completion using Mean Absolute Percentage Error (MAPE) and Root Mean Square Error (RMSE):
\begin{equation*}
\begin{aligned}
    \operatorname{MAPE} &= 100 \times \frac{1}{n} \sum_{i = 1}^{n} \frac{|y_i - \hat{y_i}|}{|y_i|}, \\
    \operatorname{RMSE} &= \sqrt{\frac{1}{n} \sum_{i = 1}^{n} (y_i - \hat{y_i})^2}, 
\end{aligned}
\end{equation*}
where $y_i$ and $\hat{y_i}$ represent actual and estimated values, and $n$ is the total number of estimated values. Lower MAPE and RMSE values reflect superior completion performance. 

To ensure fair evaluation across all compared methods, we tune hyperparameters for each method and dataset type at fixed sampling rates using a consistent strategy. For each dataset type, we select a small, representative subset of samples (e.g., \textit{bus}, \textit{mobile}, \textit{akiyo} for color videos from the YUV dataset) and perform a grid search over hyperparameter ranges recommended by the respective authors to identify optimal values. These optimal hyperparameters are then applied to all samples of the dataset type for each method. For the FBGND method~\cite{liang2024fixed}, originally designed for hyperspectral image denoising, we modified its regularizer to suit tensor completion, following~\cite{10078018}. Experiments were conducted on a machine with an Intel Xeon E5-2698 v4 CPU (2.20 GHz), 256 GB RAM, and an NVIDIA GeForce RTX 3090 (24 GB) with driver version 570.153.02.

\begin{figure*}[htbp!]
    \centering
    \includegraphics[width=.9\linewidth]{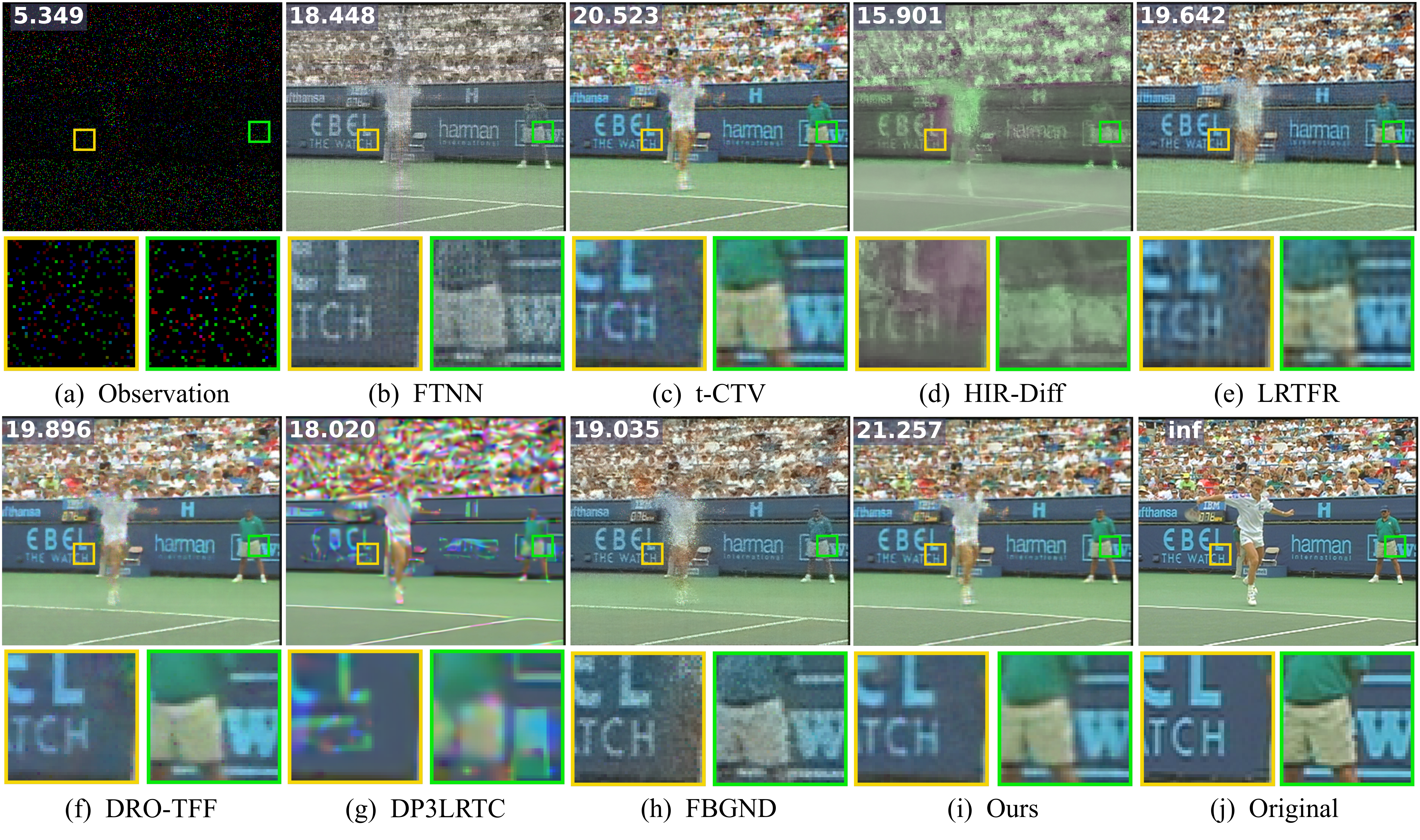}
    \caption{Completion results for the color video \textit{stefan} (SR = 0.05), showing the 25th frame of each completed video. The MPSNR is indicated in the upper-left corner of each image.}
    \label{fig:video_res}
\end{figure*}

\subsection{Implementation details} \label{subsec:implement}
We implement our method using Python with PyTorch 2.7.1 and CUDA 12.8. The invertible linear transform $\mathfrak{L}$ is the DCT, chosen for its superior empirical performance compared to the DFT~\cite{10078018}. We set a convergence threshold $\varepsilon = 10^{-4}$ and a maximum iteration limit $N = 200$. 

For the $k$-SPC deep denoiser $\operatorname{D}_\sigma$ in $\operatorname{C} = \alpha (\operatorname{Id} - \operatorname{D}_\sigma)$, we adopt the DRUNet architecture~\cite{zhang2021plug}, following~\cite{AHQS}, with pretrained weights for grayscale and color images ($k = 0.9$) from the GitHub repository\footnote{https://github.com/FizzzFizzz/Learning-Pseudo-Contractive-Denoisers-for-Inverse-Problems}. To simplify tuning and accelerate convergence, we fix the stepsize $\tau = 1$, initial penalty parameter $\rho_0 = 10^{-4}$, speed factor $\nu = 1.02$, and schedule the relaxation parameter $\lambda_t = 1$ for $t < 100$ and $\lambda_t = \frac{100}{t}$ for $t \geq 100$, satisfying the conditions $\sum \lambda_t = \infty$ and $\sum \lambda_t^2 < \infty$ in Corollary~\ref{coro:convergence}.

For MDI completion, we use the convex function $f(x) = |x|$, and denote GTCTV-DPC as Abs-TCTV-DPC in this case, with $\Gamma = \{1, 2, 4\}$, and the hyperparameters are set as follows:
\begin{itemize}
    \item \textbf{MSIs}: $N_{\text{in}} = 8$, $\sigma_0 = 0.05$, $\alpha = 1.00$.
    \item \textbf{Color videos}: $N_{\text{in}} = 5$, $\sigma_0 = 0.30$, $\alpha = 0.50$.
\end{itemize}

For spatio-temporal traffic data completion, we use the following SCAD penalty~\cite{fanvariable} within the GTCTV prior, denoted as SCAD-TCTV-DPC:
\begin{equation*}
f_{\varphi, \omega}(x) = \begin{cases}
    \varphi x, & 0 \leq x < \varphi, \\
    \frac{-x^2 + 2 \omega \varphi x - \varphi^2}{2(\omega - 1)}, & \varphi \leq x < \omega \varphi, \\
    \frac{\omega + 1}{2} \varphi^2, & x \geq \omega \varphi,
\end{cases}
\end{equation*}
with $\varphi > 0$ and $\omega > 1$. Since $f_{\varphi, \omega}$ is $\frac{1}{\omega-1}$-weakly convex, the GTCTV prior $\|\cdot\|_{\text{GTCTV}}$ is $\frac{4}{\omega-1}$-weakly convex on $\mathbb{V}$, by Lemma~\ref{lem:mu0-gtctv}. The traffic data are reshaped into a tensor of dimensions (locations/sensors $\times$ time intervals $\times \, 1 \, \times$ days), with $N_{\text{in}} = 5$, $\Gamma = \{1, 2, 4\}$ and other hyperparameters set as follows:
\begin{itemize}
    \item \textbf{Guangzhou}: $\varphi = 5.00$, $\omega = 2000$, $\sigma_0 = 0.85$, $\alpha = 1.50$.
    \item \textbf{Seattle}: $\varphi = 3.00$, $\omega = 3000$, $\sigma_0 = 0.95$, $\alpha = 2.00$.
    \item \textbf{PeMS}: $\varphi = 3.00$, $\omega = 200$, $\sigma_0 = 0.65$, $\alpha = 2.00$.
\end{itemize}

\begin{table}
    \centering
    \caption{Average quantitative results of different methods for MDI completion, with the best and second-best results highlighted in \textbf{bold} and \underline{underlined}, respectively.} \label{table:average_mdi_res}
    {\fontsize{6pt}{7pt}\selectfont
    \begin{tabular}{c cc cc cc}
        \specialrule{.1em}{.05em}{.05em}
        \specialrule{.1em}{.05em}{.05em}
        Sampling Rate & \multicolumn{2}{c}{0.05} & \multicolumn{2}{c}{0.10} & \multicolumn{2}{c}{0.20}  \\
        \specialrule{.05em}{.25em}{.25em}
        Method & MPSNR & MSSIM & MPSNR & MSSIM & MPSNR & MSSIM \\
        
        \specialrule{.05em}{.25em}{.25em}
        \multicolumn{7}{c}{Multi-Spectral Images $(256 \times 256 \times 31)$} \\
        \specialrule{.05em}{.25em}{.25em}
        
        Observation & 15.231 & 0.238 & 15.466 & 0.271 & 15.978 & 0.330 \\
        FTNN & 34.162 & 0.912 & 38.770 & 0.961 & 43.285 & 0.981 \\
        t-CTV & \underline{37.839} & \underline{0.960} & \underline{41.430} & \underline{0.978} & \underline{45.549} & \underline{0.990} \\
        HIR-Diff & 23.031 & 0.736 & 24.738 & 0.769 & 26.681 & 0.829 \\
        LRTFR & 37.232 & 0.951 & 40.519 & 0.971 & 42.614 & 0.972 \\
        DRO-TFF & 37.462 & 0.964 & 40.407 & 0.975 & 43.685 & 0.988 \\
        DP3LRTC & 33.886 & 0.942 & 37.127 & 0.966 & 40.355 & 0.980 \\
        FBGND & 27.546 & 0.734 & 30.195 & 0.794 & 32.400 & 0.852 \\
        \rowcolor{Platinum}
        Ours & \textbf{38.556} & \textbf{0.970} & \textbf{41.856} & \textbf{0.982} & \textbf{45.748} & \textbf{0.990} \\
        
        \specialrule{.05em}{.25em}{.25em}
        \multicolumn{7}{c}{Color Videos $(288 \times 352 \times 3 \times 50)$} \\
        \specialrule{.05em}{.25em}{.25em}
        
        Observation & 6.457 & 0.021 & 6.692 & 0.034 & 7.203 & 0.058 \\
        FTNN & 22.981 & 0.699 & 25.611 & 0.810 & 28.560 & 0.892 \\
        t-CTV & \underline{26.857} & \underline{0.766} & \underline{29.099} & \underline{0.840} & \underline{32.368} & \underline{0.914} \\
        HIR-Diff & 19.733 & 0.558 & 20.513 & 0.577 & 21.159 & 0.610 \\
        LRTFR & 24.856 & 0.681 & 26.179 & 0.730 & 27.300 & 0.767 \\
        DRO-TFF & 25.394 & 0.709 & 27.001 & 0.785 & 28.969 & 0.846 \\
        DP3LRTC & 23.696 & 0.780 & 25.907 & 0.853 & 28.609 & 0.915 \\
        FBGND & 24.144 & 0.733 & 26.325 & 0.811 & 27.799 & 0.854 \\
        \rowcolor{Platinum}
        Ours & \textbf{27.506} & \textbf{0.809} & \textbf{29.694} & \textbf{0.873} & \textbf{32.813} & \textbf{0.931} \\
        
        \specialrule{.1em}{.05em}{.05em}
        \specialrule{.1em}{.05em}{.05em}
    \end{tabular}}
\end{table}

\subsection{Completion performence} \label{subsec:completion}
In Table~\ref{table:average_mdi_res}, we present the average quantitative results for MDI completion across various methods. Table~\ref{table:average_mdi_res} shows that Abs-TCTV-DPC performs particularly well on more intricate data, such as MSIs and color videos. Notably, at a low sampling rate of 0.05, Abs-TCTV-DPC achieves an average MPSNR improvement of 0.717 dB for MSIs and 0.649 dB for color videos over the second-best method, highlighting its superior completion performance.

\begin{figure*}[htbp]
    \centering
    \includegraphics[width=.9\linewidth]{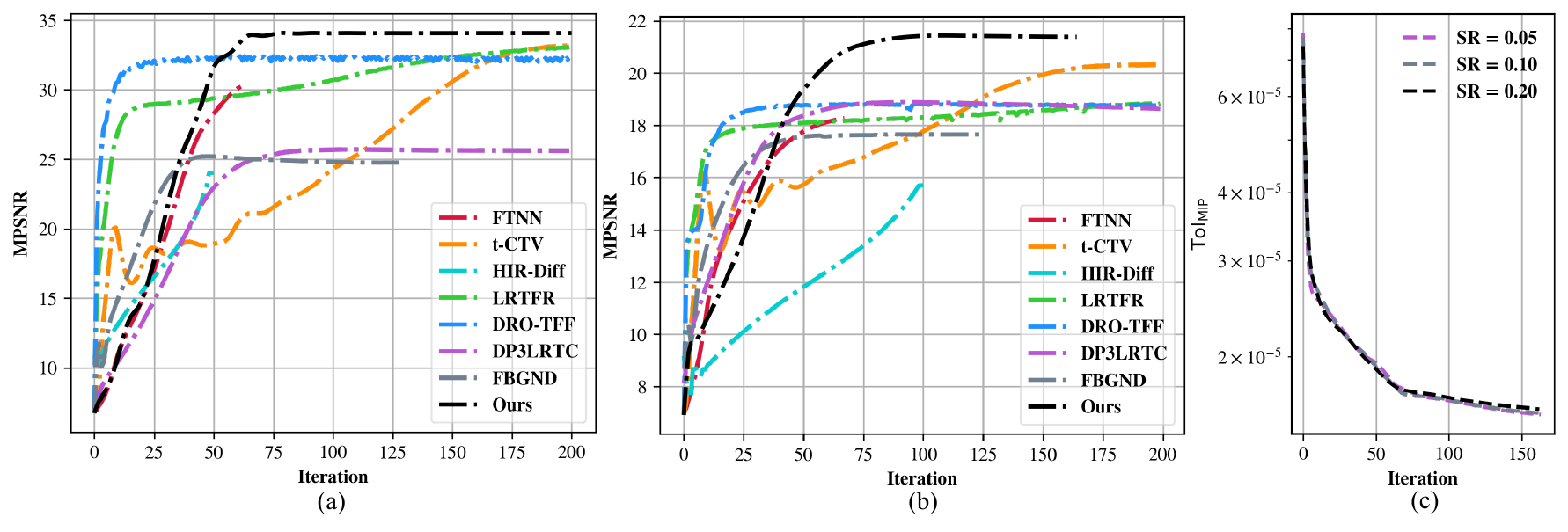}
    \caption{(a) and (b): MPSNR curves with the $x$-axis denoting the iteration number for the MSI \textit{paints} (a) and the video \textit{bus} (b) at SR = 0.05. Note that for LRTFR and DRO-TFF, which run up to 3000 iterations, we record the MPSNR every 15 iterations. (c): Convergence performance within the monotone inclusion paradigm for the video \textit{bus} at SR = 0.05, 0.10, and 0.20.}
     \label{fig:psnr_curve}
\end{figure*}

In Fig.~\ref{fig:video_res}, we present the visual results for the color video \textit{stefan} at SR = 0.05. The visual results indicate that Abs-TCTV-DPC excels in preserving both the overall structural coherence and intricate local patterns, delivering clearer and more detailed reconstructions even at low sampling rates. Additionally, Fig.~\ref{fig:psnr_curve} illustrates the MPSNR curves for the MSI \textit{paints} and the video \textit{bus} at SR = 0.05 across all competing methods. The curves illustrate that GTCTV-DPC achieves stable and progressive MPSNR improvements over iterations, highlighting its robustness and consistency.

To assess convergence from the monotone-inclusion perspective, Fig.~\ref{fig:psnr_curve}(c) reports the following residual
\begin{equation*}
    \operatorname{Tol}_{\mathrm{MIP}}(t)
    \;:=\;
    \frac{\big\| \big(\operatorname{A} + \operatorname{B} + \operatorname{C}\big)\big(\mathcal{X}_t^{\operatorname{A}}\big) \big\|_F}
         {\prod_{i=1}^N n_i},
\end{equation*}
where $\operatorname{A}$, $\operatorname{B}$, and $\operatorname{C}$ correspond to the operators in Eq.~\eqref{eq:concrete_set}, using the color video \textit{bus} at SR = 0.05, 0.10, and 0.20 as examples. Here $\operatorname{A}=\partial\delta_{\mathcal{Y},\Omega}$ is the subdifferential of the data-consistency indicator; by Example~3.5 in \cite{beck2017first} we have
\begin{equation*}
\begin{aligned}
&\operatorname{A} \mathcal{X}_t^{\operatorname{A}} = \partial \delta_{\mathcal{Y}, \Omega} \left(\mathcal{X}_t^{\operatorname{A}}\right) \\
=& \left\{ \mathcal{X} \in \mathbb{V} \mid \langle \mathcal{X}, \mathcal{Z} -\mathcal{X}_t^{\operatorname{A}} \rangle \le 0, \right.\\
&\left.\forall \mathcal{Z} \in \mathbb{V} \text{, such that } \mathcal{P}_{\Omega}(\mathcal{Z}) = \mathcal{P}_{\Omega}(\mathcal{Y}) \right\}.
\end{aligned}
\end{equation*}
We select $\mathcal{O} \in \operatorname{A} \mathcal{X}_t^{\operatorname{A}}$ as the subgradient satisfying this inclusion. For the GTCTV term, $\operatorname{B} \mathcal{X}_t^{\operatorname{A}}$ is obtained by evaluating a representative subgradient of the GTCTV penalty at $\mathcal{X}_t^{\operatorname{A}}$. In practice, we compute $\left\| \mathcal{X}_t^{\operatorname{A}} \right\|_{\text{GTCTV}}$ and employ the `torch.autograd` module to obtain its gradient automatically. Finally, $\operatorname{C}(\mathcal{X})=\alpha\big(\operatorname{Id}-\operatorname{D}_\sigma\big)(\mathcal{X})$ is evaluated directly using the DPC denoiser. Fig.~\ref{fig:psnr_curve}(c) demonstrates that $\operatorname{Tol}_{\mathrm{MIP}}(t)$ decays to near zero as iterations progress, indicating convergence of the iterates to a solution of the monotone inclusion problem in Eq.~\eqref{eq:Proposed-MIP}. Notably, these experiments incorporate the practical acceleration heuristics outlined in Algorithm~\ref{alg:main} and section~\ref{subsec:implement}. Despite these speedups, the monotone inclusion residual converges, providing empirical support for the theoretical global convergence result.

\begin{table}[htbp]
    \centering
    \caption{The quantitative results of different methods for traffic data completion.} \label{table:traffic_res}
    {\fontsize{6pt}{7pt}\selectfont
    \begin{tabular}{c cc cc cc}
        \specialrule{.1em}{.05em}{.05em}
        \specialrule{.1em}{.05em}{.05em}
        Sampling Rate & \multicolumn{2}{c}{0.30} & \multicolumn{2}{c}{0.50} & \multicolumn{2}{c}{0.70}  \\
        \specialrule{.05em}{.25em}{.25em}
        Method & MAPE & RMSE & MAPE & RMSE & MAPE & RMSE \\
        
        \specialrule{.05em}{.25em}{.25em}
        \multicolumn{7}{c}{Guangzhou $(214 \times 144 \times 61)$} \\
        \specialrule{.05em}{.25em}{.25em}

        BATF & 8.55 & 3.70 & 8.36 & 3.61 & 8.30 & 3.59 \\
        BTMF & 8.65 & 3.69 & 7.89 & 3.40 & 7.44 & 3.22 \\
        LRTC-TNN & 8.39 & 3.60 & 7.66 & 3.29 & 7.02 & 3.02 \\
        LRTC-TSpN & 8.62 & 3.66 & 7.77 & 3.31 & 7.06 & 3.01 \\
        LSTC-Tubal & {8.21} & \underline{3.47} & \underline{7.26} & \underline{3.10} & \underline{6.64} & \underline{2.85} \\
        LATC & 8.43 & 3.62 & 7.68 & 3.29 & 7.04 & 3.02 \\
        LRTFR & \underline{8.13} & 3.51 & 7.23 & 3.14 & 6.70 & 2.90 \\
        \rowcolor{Platinum}
        Ours & \textbf{7.95} & \textbf{3.47} & \textbf{6.98} & \textbf{3.07} & \textbf{6.38} & \textbf{2.81} \\
        
        \specialrule{.05em}{.25em}{.25em}
        \multicolumn{7}{c}{Seattle $(323 \times 288 \times 28)$} \\
        \specialrule{.05em}{.25em}{.25em}
        
        BATF & 7.38 & 4.46 & 7.20 & 4.37 & 7.16 & 4.34 \\
        BTMF & 6.22 & 3.86 & 5.80 & 3.65 & 5.64 & 3.57 \\
        LRTC-TNN & 6.56 & 3.96 & 5.56 & 3.47 & 4.95 & 3.16 \\
        LRTC-TSpN & 6.44 & 3.93 & 5.53 & 3.47 & \textbf{4.75} & \textbf{3.08} \\
        LSTC-Tubal & 6.93 & 4.09 & 6.12 & 3.70 & 5.65 & 3.48 \\
        LATC & \underline{6.06} & \underline{3.77} & \underline{5.34} & \underline{3.39} & {4.90} & {3.15} \\
        LRTFR & 7.20 & 4.42 & 6.57 & 4.01 & 6.21 & 3.80 \\
        \rowcolor{Platinum}
        Ours & \textbf{5.95} & \textbf{3.73} & \textbf{5.28} & \textbf{3.37} & \underline{4.87} & \underline{3.14} \\
        
        \specialrule{.05em}{.25em}{.25em}
        \multicolumn{7}{c}{PeMS $(228 \times 288 \times 44)$} \\
        \specialrule{.05em}{.25em}{.25em}
        
        BATF & 6.96 & 4.78 & 6.83 & 4.71 & 6.82 & 4.68 \\
        BTMF & 5.41 & 3.87 & 4.92 & 3.62 & 4.64 & 3.50 \\
        LRTC-TNN & 5.94 & 4.15 & 4.45 & 3.13 & 3.45 & 2.44 \\
        LRTC-TSpN & 4.63 & \underline{3.21} & \underline{3.42} & \textbf{2.42} & \textbf{2.74} & \textbf{1.96} \\
        LSTC-Tubal & \underline{4.56} & 3.22 & 3.51 & 2.51 & 2.95 & 2.09 \\
        LATC & 5.10 & 3.58 & 4.00 & 2.84 & 3.30 & 2.35 \\
        LRTFR & 6.33 & 4.05 & 4.93 & 3.25 & 4.17 & 2.72 \\
        \rowcolor{Platinum}
        Ours & \textbf{4.29} & \textbf{3.13} & \textbf{3.38} & \underline{2.44} & \underline{2.88} & \underline{2.05} \\
        
        \specialrule{.1em}{.05em}{.05em}
        \specialrule{.1em}{.05em}{.05em}
    \end{tabular}}
\end{table}

\begin{figure}[htbp!]
    \centering
    \includegraphics[width=1.\linewidth]{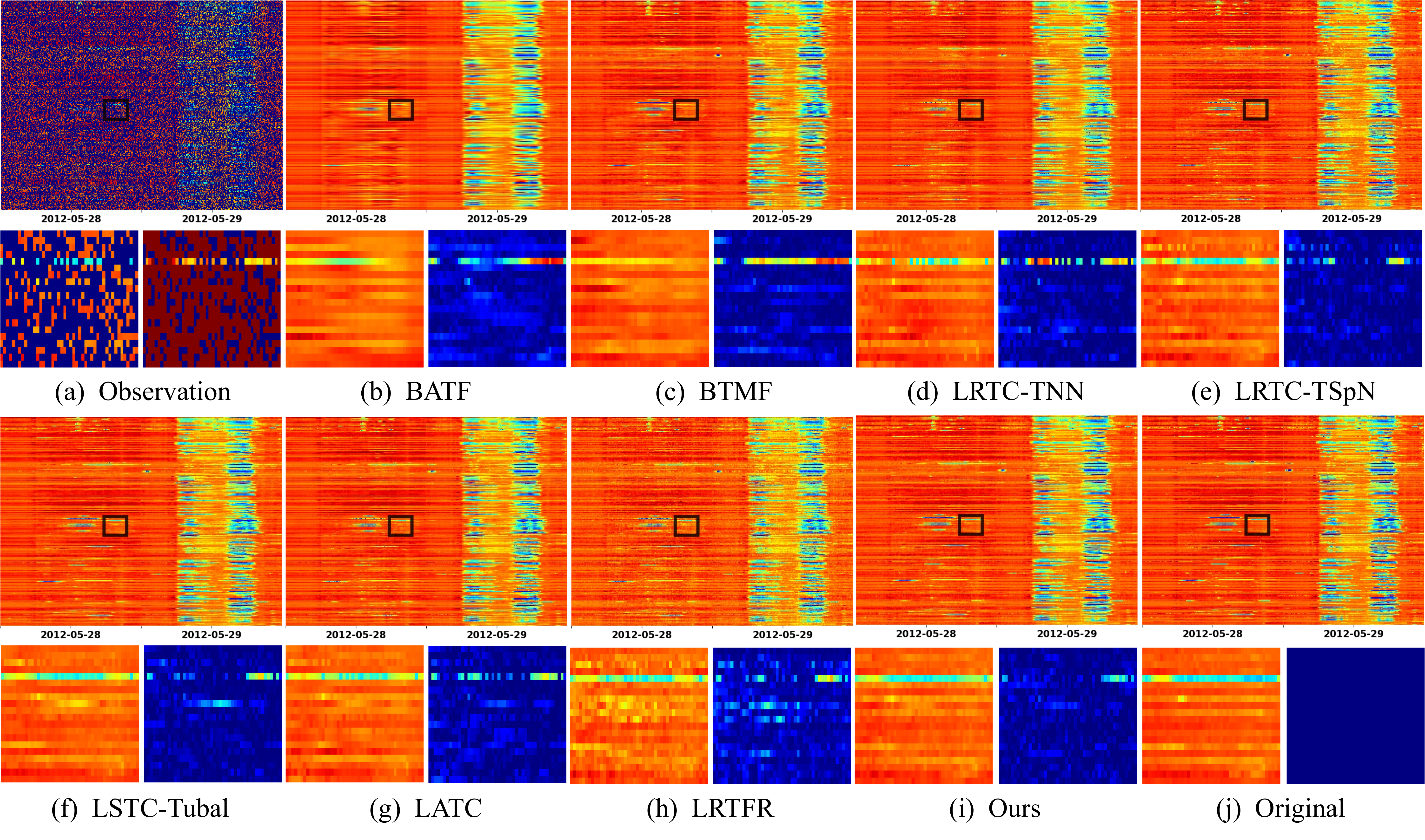}
    \caption{Completion results for the traffic data \textit{PeMS} (SR = 0.30), showing the enlarged views of the box regions alongside corresponding residual components.}
    \label{fig:traffic_res}
\end{figure}

Table~\ref{table:traffic_res} summarizes the quantitative results across multiple methods. The results indicate that SCAD-TCTV-DPC achieves competitive imputation performance compared to both low-rank methods with nonconvex penalties and multi-prior approaches, particularly at a low sampling rate of 0.30. Fig.~\ref{fig:traffic_res} presents the imputation results for two selected days from PeMS at SR = 0.30, demonstrating the superior performance of our method in effectively preserving the global structure while maintaining intricate local details of the original traffic data. In section~S.\Rmnum{2}.B of the supplement, we further evaluate our method on color image completion. Additional visual results are provided in section~S.\Rmnum{2}.C.

\subsection{Discussions}
\subsubsection{Ablation studies} \label{subsec:ablation_studies}
We perform two ablation studies: (i) to evaluate the contributions of the GTCTV prior and the DPC denoiser in GTCTV-DPC for MDI completion; and (ii) to assess the effect of the SCAD penalty for spatio-temporal traffic-data completion. For the MDI completion ablation, we consider three variants of GTCTV-DPC:
\begin{itemize}
    \item \textbf{TNN-DPC}: Replace the GTCTV prior with TNN to test the effectiveness of GTCTV, while retaining the DPC denoiser; 
    \item \textbf{Abs-TCTV-DNE}: GTCTV with $f(x)=|x|$ paired with a deep non-expansive (Definition~\ref{def:ne}, NE) denoiser to test the effect of the denoiser assumption; 
    \item \textbf{Abs-TCTV-DFNE}: GTCTV with $f(x)=|x|$ and a deep firm non-expansive (Definition~\ref{def:firm}, FNE) denoiser~\cite{pesquet2021learning}.
\end{itemize}

These variants are tested on five YUV color videos (\textit{bus}, \textit{mobile}, \textit{akiyo}, \textit{mother-daughter}, \textit{tempete}) at sampling rates of 0.05, 0.10, and 0.20. Hyperparameters follow GTCTV-DPC for color videos ($N_{\text{in}} = 5$, $\Gamma = \{1, 2, 4\}$, $\sigma_0 = 0.30$, $\alpha = 0.50$), except for TNN-DPC, which uses an initial TNN threshold of $1$ and divides the speed factor $\nu = 1.02$ by the iteration number. All denoisers (DPC, NE, FNE) use the DRUNet architecture~\cite{zhang2021plug}, with pretrained weights for NE\footnote{https://github.com/FizzzFizzz/New-baseline-for-DRUNet-under-different-assumptions} and FNE\footnote{https://github.com/basp-group/PnP-MMO-imaging}. Completion performance is evaluated using MPSNR and MSSIM, with results averaged over five runs reported in Table~\ref{table:average_ablation_res}.
\begin{table}
    \centering
    \caption{Average MPSNR and MSSIM for GTCTV-DPC variants on five YUV color videos at sampling rates 0.05, 0.10, and 0.20.} \label{table:average_ablation_res}
    {\fontsize{6pt}{7pt}\selectfont
    \begin{tabular}{c cc cc cc}
        \specialrule{.1em}{.05em}{.05em}
        \specialrule{.1em}{.05em}{.05em}
        Sampling Rate & \multicolumn{2}{c}{0.05} & \multicolumn{2}{c}{0.10} & \multicolumn{2}{c}{0.20}  \\
        \specialrule{.05em}{.25em}{.25em}
        Method & MPSNR & MSSIM & MPSNR & MSSIM & MPSNR & MSSIM \\

        \specialrule{.05em}{.25em}{.25em}
        
        Observation & 7.503 & 0.027 & 7.738 & 0.043 & 8.249 & 0.078 \\
        TNN-DPC & 23.480 & 0.633 & 26.985 & 0.774 & 30.651 & 0.883 \\
        Abs-TCTV-DNE & \underline{26.869} & \underline{0.781} & \underline{29.365} & \underline{0.856} & \underline{32.780} & \underline{0.924} \\
        Abs-TCTV-DFNE & 26.442 & 0.765 & 28.943 & 0.843 & 32.341 & 0.917 \\
        \rowcolor{Platinum}
        Abs-TCTV-DPC & \textbf{27.274} & \textbf{0.796} & \textbf{29.769} & \textbf{0.868} & \textbf{33.058} & \textbf{0.931} \\
        
        \specialrule{.1em}{.05em}{.05em}
        \specialrule{.1em}{.05em}{.05em}
    \end{tabular}}
\end{table}

\begin{figure}
    \centering
    \includegraphics[width=1.\linewidth]{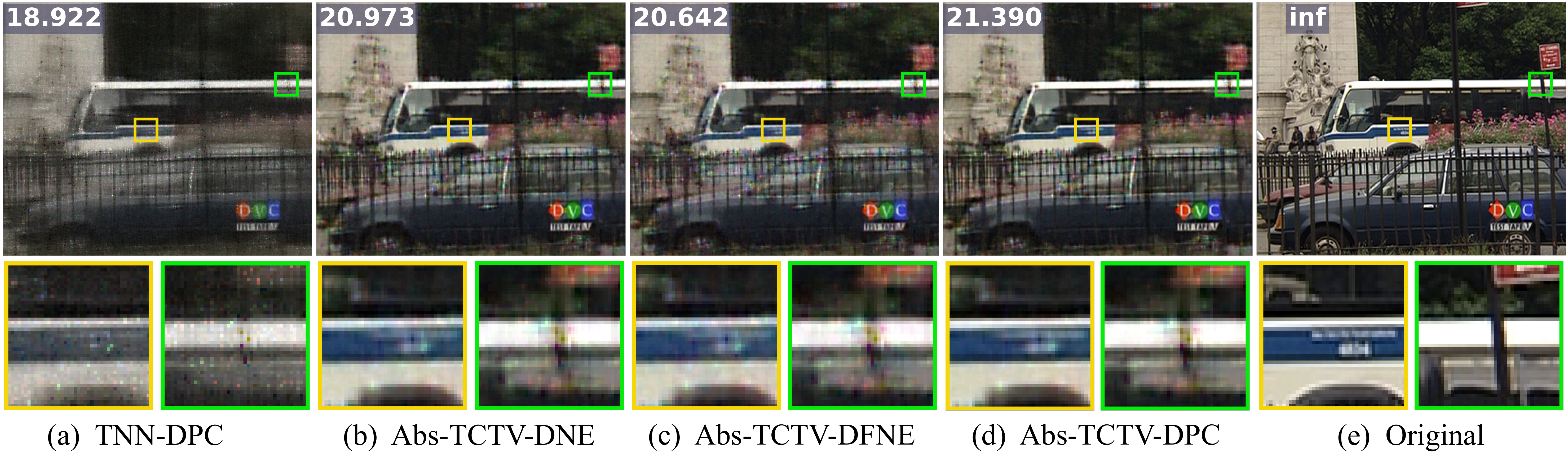}
    \caption{The results of color video completion on \textit{bus} (SR = 0.05) by different variants of the proposed method.}
    \label{fig:ablation_video_res}
\end{figure}

Table~\ref{table:average_ablation_res} shows that Abs-TCTV-DPC outperforms other variants across all sampling rates, with Abs-TCTV-DNE achieving the second-best performance. Fig.~\ref{fig:ablation_video_res} shows the visual comparisons for the \textit{bus} video (SR = 0.05). The GTCTV prior, incorporating the spatial and temporal gradient information, captures the holistic tensor structure more effectively than TNN, leading to higher MPSNR and MSSIM. Furthermore, as discussed in section~\ref{subsec:functional_analysis}, the hierarchy of denoiser assumptions (FNE $\implies$ NE $\implies$ PC, Eq.~\eqref{eq:relation_operators}) indicates that the less restrictive PC assumption in GTCTV-DPC enables better completion performance compared to the stronger NE and FNE assumptions, as evidenced by the results. 

\begin{table}[htbp]
    \centering
    \caption{The quantitative results for traffic data completion.} \label{table:traffic_ablation_res}
    {\fontsize{6pt}{7pt}\selectfont
    \begin{tabular}{c cc cc cc}
        \specialrule{.1em}{.05em}{.05em}
        \specialrule{.1em}{.05em}{.05em}
        Sampling Rate & \multicolumn{2}{c}{0.30} & \multicolumn{2}{c}{0.50} & \multicolumn{2}{c}{0.70}  \\
        \specialrule{.05em}{.25em}{.25em}
        Method & MAPE & RMSE & MAPE & RMSE & MAPE & RMSE \\
        
        \specialrule{.05em}{.25em}{.25em}
        \multicolumn{7}{c}{Guangzhou $(214 \times 144 \times 61)$} \\
        \specialrule{.05em}{.25em}{.25em}

        Abs-TCTV-DPC & 9.25 & 4.02 & 7.29 & 3.20 & 6.38 & 2.82 \\
        SCAD-TCTV-DPC & \textbf{7.95} & \textbf{3.47} & \textbf{6.98} & \textbf{3.07} & \textbf{6.38} & \textbf{2.81} \\
        
        \specialrule{.05em}{.25em}{.25em}
        \multicolumn{7}{c}{Seattle $(323 \times 288 \times 28)$} \\
        \specialrule{.05em}{.25em}{.25em}

        Abs-TCTV-DPC & 7.43 & 4.68 & 6.59 & 4.27 & 6.15 & 4.03 \\
        SCAD-TCTV-DPC & \textbf{5.95} & \textbf{3.73} & \textbf{5.28} & \textbf{3.37} & \textbf{4.87} & \textbf{3.14} \\
        
        \specialrule{.05em}{.25em}{.25em}
        \multicolumn{7}{c}{PeMS $(228 \times 288 \times 44)$} \\
        \specialrule{.05em}{.25em}{.25em}

        Abs-TCTV-DPC & 5.02 & 3.65 & 3.58 & 2.62 & 2.92 & 2.10 \\
        SCAD-TCTV-DPC & \textbf{4.29} & \textbf{3.13} & \textbf{3.38} & \textbf{2.44} & \textbf{2.88} & \textbf{2.05} \\
        
        \specialrule{.1em}{.05em}{.05em}
        \specialrule{.1em}{.05em}{.05em}
    \end{tabular}}
\end{table}

To assess the significance of the SCAD penalty in GTCTV-DPC for spatio-temporal traffic data completion, we compare SCAD-TCTV-DPC against Abs-TCTV-DPC. We test these variants on the Guangzhou ($214 \times 144 \times 61$), Seattle ($323 \times 288 \times 28$), and PeMS ($228 \times 288 \times 44$) datasets at sampling rates of 0.30, 0.50, and 0.70. Completion performance is evaluated using MAPE and RMSE, with results reported in Table~\ref{table:traffic_ablation_res}. Table~\ref{table:traffic_ablation_res} shows that SCAD-TCTV-DPC outperforms Abs-TCTV-DPC across all datasets and sampling rates, with lower MAPE and RMSE. The nonconvex SCAD penalty promotes sparsity and captures complex spatio-temporal patterns in traffic data more effectively than the convex absolute value function, demonstrating its suitability for such datasets.

\begin{figure}[htbp!]
    \centering
    \includegraphics[width=.95\linewidth]{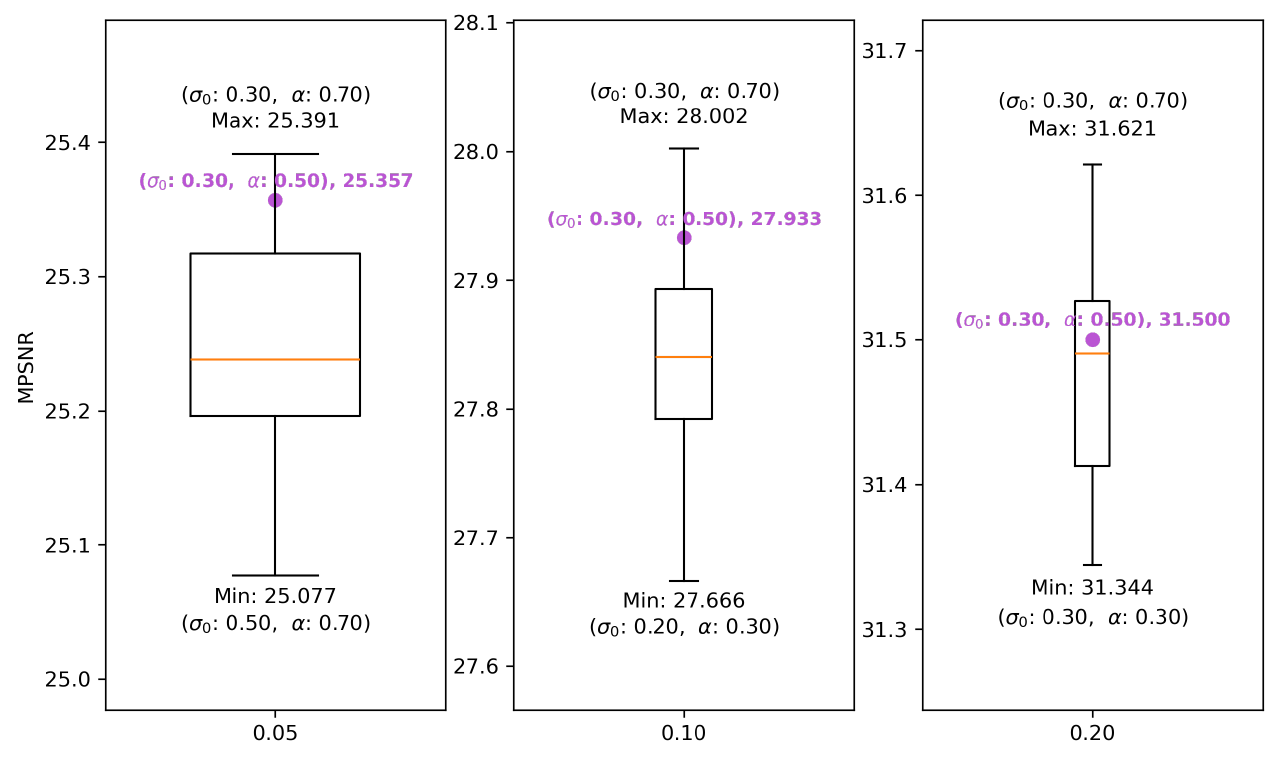}
    \caption{Boxplot of average MPSNR for each $\sigma_0$-$\alpha$ combination for Abs-TCTV-DPC on the YUV color video subset (\textit{bus}, \textit{mobile}, \textit{akiyo}) across sampling rates 0.05, 0.10, and 0.20. The selected combination ($\sigma_0 = 0.30$, $\alpha = 0.50$) is highlighted.}
    \label{fig:hyper_params}
\end{figure}

\subsubsection{Hyper-parameters tuning strategy} \label{subsec:hyper_tuning}
As introduced in section~\ref{subsec:expe_sets}, for each dataset type with a fixed sampling rate, we choose the same small samples from selected datas, and utilize the grid search to find the optimal hyper-parameters for the small samples, and finally adopt the optimal hyper-parameters for all selected datas for all compared methods. Here we take the completion of color videos as an instance to introduce the procession of hyper-parameters tunning strategy for GTCTV-DPC.

For color videos, we select \textit{bus}, \textit{mobile}, and \textit{akiyo} from the YUV dataset, as shown in Fig.~S.2 from the supplement, and perform grid search at sampling rates of 0.05, 0.10, and 0.20. For GTCTV-DPC, two hyper-parameters need to be carefully tuned: the initial denoising strength $\sigma_0$ and the weighting factor $\alpha$ in the DPC denoiser $\operatorname{C} = \alpha (\operatorname{Id} - \operatorname{D}_\sigma)$. We search over $\sigma_0 \in \{0.20, 0.30, 0.50, 0.70\}$ and $\alpha \in \{0.30, 0.50, 0.70\}$, resulting in 12 combinations. Fig.~\ref{fig:hyper_params} presents a boxplot of the average MPSNR for each combination of $\sigma_0$ and $\alpha$ across all sampling rates, illustrating the robustness of GTCTV-DPC's completion performance to different hyper-parameter choices. For color videos, we select $\sigma_0 = 0.30$ and $\alpha = 0.50$, as these values balance high MPSNR and MSSIM with computational efficiency. These selected hyper-parameters were then applied to all color video datasets in our experiments.

\begin{table}[htbp]
    \centering
    \caption{Computational time (mean $\pm$ std, in seconds).} \label{table:used_time}
    \footnotesize{
    \begin{tabular}{c | c c c}
        \specialrule{.1em}{.05em}{.05em}
        \specialrule{.1em}{.05em}{.05em}
        Method & \multicolumn{3}{c}{Used Time (s)} \\
        \specialrule{.05em}{.25em}{.25em}
        HIR-Diff & 63.25 & $\pm$ & 8.42 \\
        LRTFR & 20.84 & $\pm$ & 0.31 \\
        DRO-TFF & 21.65 & $\pm$ & 0.47 \\
        \specialrule{.05em}{.25em}{.25em}
        FTNN & 10994.82 & $\pm$ & 107.80 \\
        t-CTV & 1540.58 & $\pm$ & 5.85 \\
        DP3LRTC & 1207.08 & $\pm$ & 11.13 \\
        FBGND & 5543.76 & $\pm$ & 9.29 \\
        \rowcolor{Platinum}
        Ours & 2125.01 & $\pm$ & 4.05 \\
        \specialrule{.1em}{.05em}{.05em}
        \specialrule{.1em}{.05em}{.05em}
    \end{tabular}}
\end{table}

\subsubsection{Computational performance} \label{subsec:computational}
\begin{figure}[htbp!]
    \centering
    \includegraphics[width=.8\linewidth]{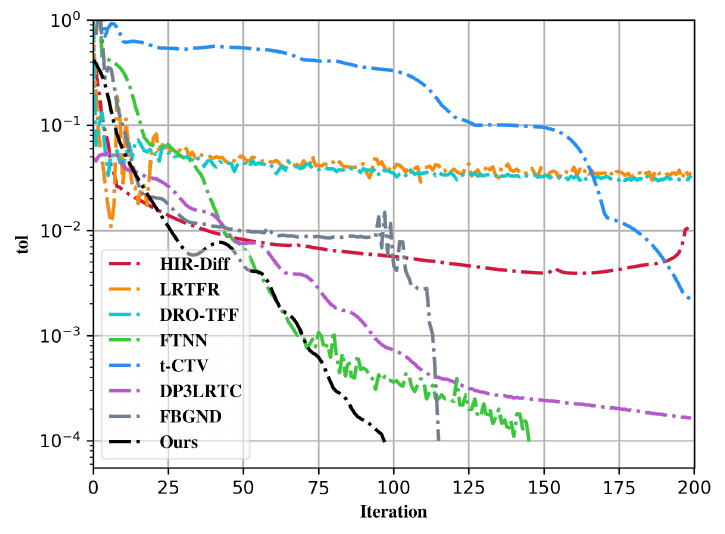}
    \caption{Residual error (tolerance) versus iteration number for compared methods on the \textit{akiyo} video ($288 \times 352 \times 3 \times 50$, sampling rate 0.20). For LRTFR and DRO-TFF, tolerance is recorded every 15 iterations due to their 3000-iteration limit.}
    \label{fig:tol_curve}
\end{figure}

To test computational performance, we use the \textit{akiyo} color video ($288 \times 352 \times 3 \times 50$) from the YUV dataset with a sampling rate of 0.20. Each method is run independently five times, and Table~\ref{table:used_time} reports the mean and standard deviation of computational time (in seconds). Fig.~\ref{fig:tol_curve} shows the residual error (tolerance) as a function of iteration number for each method. For LRTFR and DRO-TFF, which are run for up to 3000 iterations, the tolerance is recorded every 15 iterations.

Self-supervised deep-learning methods (HIR-Diff, LRTFR, DRO-TFF) leverage lightweight architectures and GPU acceleration, yielding high computational efficiency. However, their completion performance, as measured by MPSNR and MSSIM, remains inferior to that of GTCTV-DPC and t-CTV. The proposed GTCTV-DPC method, which integrates the GTCTV prior, exhibits rapid convergence in terms of residual error, as illustrated in Fig.~\ref{fig:tol_curve}. Nevertheless, due to the multiple t-SVD computations required per iteration, it incurs higher computational costs. Despite this, GTCTV-DPC achieves consistently superior completion performance, as shown in section~\ref{subsec:completion}. Moreover, at a lower sampling rate of 0.05, GTCTV-DPC demonstrates further improvement in completion quality, indicating its strong adaptability under highly undersampled conditions. In the future, we will explore the use of randomized SVD techniques to accelerate t-SVD computations and enhance scalability.

\section{Conclusions} \label{sec:conclusion}
In this work, we present a novel tensor completion method within the monotone inclusion paradigm. To effectively capture global structure, we generalize the t-CTV prior with a weakly convex penalty and rigorously established its weak convexity in Lemma~\ref{lem:mu0-gtctv}. To preserve intricate local details, we incorporate DPC denoisers~\cite{AHQS} and establish their connection with $\beta$-cocoercivity in Lemma~\ref{lem:spc and coco}. Leveraging the DYS scheme, we derive the GTCTV-DPC method in Algorithm~\ref{alg:main}. 

A key theoretical contribution of this work lies in the convergence analysis of the GTCTV-DPC method. Starting from the general MIP in Eq.~\eqref{eq:MIP}, we showed that the associated operator $\operatorname{T}$ in Eq.~\eqref{eq:T in DYS} is SPC and extended the admissible stepsize range, along with explicit conditions for relaxation parameters. This yields Corollary~\ref{coro:convergence}, which establishes the global convergence of Algorithm~\ref{alg:main} to a solution of the proposed model in Eq.~\eqref{eq:Proposed-MIP}. 

Empirical results on MDI and traffic datasets further demonstrate the superior performance and strong visual fidelity of the proposed method. For instance, at a sampling rate of 0.05 for MDI completion, GTCTV-DPC achieves an average MPSNR that surpasses the second-best method by 0.717 dB, and 0.649 dB for MSIs, and color videos, respectively. 

Despite these advances, the reliance on multiple t-SVD operations results in high computational complexity, particularly for large-scale tensors~\cite{HOU2025128912}, as discussed in section~\ref{subsec:computational}. Future research will investigate randomized SVD algorithms~\cite{doi:10.1137/090771806} to enhance computational efficiency and investigate the adaptability of our framework to a wider range of tensor recovery tasks with tailored weakly convex penalties to further enhance practical applicability.


\bibliographystyle{IEEEtran}

\bibliography{references}

\end{document}